\magnification=\magstep1
\tolerance=10000
\hbadness=10000
\vbadness=10000

\def\Hom{{\rm Hom}}
\def\Z{{\bf Z}}   
\def\Q{{\bf Q}}   
\def\C{{\bf C}}   

\proclaim Vertex algebras. \hfill 10 June 1997, corrected 19 June, 
26 June, 26 September.

Richard E. Borcherds,
\footnote{$^*$}{Supported by NSF grant
DMS-9401186 and  a Royal Society research professorship.}

D.P.M.M.S.,
16 Mill Lane,
Cambridge,
CB2 1SB,
England.

\bigskip

\beginsection Contents.

\item{1.} Introduction.
\item{2.} Formal groups and Hopf algebras.
\item{3.} Vertex groups
\item{4.} Relaxed multilinear categories
\item{5.} The representations of a vertex group.
\item{6.} $G$ vertex algebras
\item{7.} Vertex operators
\item{8.} Free field $G$ vertex algebras.
\item{9.} The main identity
\item{10.} The Yang-Baxter equation
\item{11.} Cohomology of $G$ vertex algebras.
\item{12.} Relation to quantum field theory.

\proclaim 1.~Introduction.

In this paper we try to define the higher dimensional analogues of
vertex algebras.  In other words we define algebras which we hope have
the same relation to higher dimensional quantum field theories that
vertex algebras have to one dimensional quantum field theories (or to
``chiral halves'' of two dimensional conformal field theories).

The ideas in this paper are not really in final form. On the other
hand, I have been rewriting versions of this paper for about 12 years, 
and it has by now become obvious that it  never will attain a satisfactory
final form. I apologize in advance for the resulting chaotic and
incomplete nature
of some of the sections.

The main ideas of this paper are as follows. We first define vertex
groups, which can be thought of roughly as groups together with some
allowable singularities for functions on the group. We look at the
category of representations of the group, and redefine multilinear
maps between representations to allow maps with certain sorts of
singularities.  The result of this is that the composition of
multilinear maps is no longer a multilinear map, but it is
sufficiently close to a multilinear map that we can still sensibly
define analogues of associativity, $G$-algebras, commutative
$G$-algebras, and so on.  We will refer to the analogues of $G$
algebraic structures in this relaxed multilinear category as
$G$ vertex algebraic structures; for example we can define $G$ vertex
(associative) algebras, $G$ vertex Lie algebras, and so on. We will
mostly be interested in commutative $G$ vertex algebras, which are the
analogues of commutative rings acted on by the ``group'' $G$.  We show
that vertex algebras are exactly the same as commutative associative
$G$ vertex algebras, where $G$ is the simplest nontrivial example of a
vertex group. It is now obvious how to define higher dimensional
analogues of vertex algebras: they are just the commutative
associative $G$ vertex algebras for a higher dimensional vertex group
$G$.

Many examples of commutative $G$ vertex algebras come from quantum
field theories. Very roughly, a quantum field theory ought to give a
$G$ algebra if and only if it has ``good'' operator product
expansions, and the quantum field theory is then determined by its
$G$ vertex algebra.  The underlying group of $G$ will be the group of
translations of spacetime, or some larger group such as the Poincare
group, and the ``allowable'' singularities on $G$ are the
singularities that appear in the operator product expansions.  Any
$G$ vertex algebra automatically has ``good'' operator product
expansions; in fact, $G$ vertex algebras are just a formalization of
what physicists do with operator product expansions.  A simple example
of a higher dimensional commutative $G$ vertex algebra is the space of
Wick polynomials in a free field, which is the $G$ vertex algebra of a
free quantum field.  The commutativity of the $G$ algebra is of course
very closely related to the ``locality'' property of a quantum field
theory.  On the other hand not all $G$ vertex algebras come from
quantum field theories, because there is usually no sensible way to
construct a Hilbert space from a $G$ vertex algebra.

$G$ vertex algebras are also closely related to operator product
expansion (OPE) algebras. The latter do not seem to be rigorously
defined anywhere, and $G$ vertex algebras can be thought of
as an attempt to give a
rigorous definition.

The simplest examples of $G$ vertex algebras are associative algebras
acted on by the underlying group of the vertex group $G$. In general a
$G$ vertex algebra can be thought of as a sort of associative algebra
acted on by a group $G$, except that the ring multiplication has
singularities and is not defined everywhere.  $G$ vertex algebras have
most of the formal properties of associative algebras: for example, we
can define left and right ideals, homomorphisms, tensor products,
commutative $G$ vertex algebras, and so on. We can also define left
and right modules over $G$ vertex algebras and multilinear maps
between them, and these multilinear maps are sometimes represented by
``tensor products''.  {}From this viewpoint quantization can be
described as follows: quantization means deforming an honest
commutative algebra acted on by the underlying group of $G$ into a
commutative $G$ vertex algebra over a vertex group $G$. The
commutative algebra we start with will usually be an associative
algebra generated by classical fields and their derivatives (of all orders).
The correlation functions of the quantized theory can be recovered
from the $G$ vertex algebra if we are given a trace on the $G$ vertex
algebra.

Another way to think of $G$ vertex algebras is as follows. 
If $V$ is an associative ring acted on by a group $G$,
then for any fixed  $v_i\in V$ and any $g_i\in G$
we can think of $v_1^{g_1}v_2^{g_2}\cdots v_i^{g_i}$
as being a function on $G^i$ with values in $V$, and it is not difficult
to write down a set of axioms for these functions equivalent to the
axioms for an associative algebra acted on by $G$. The definition of
a $G$ vertex algebra is now in principle very easy: we
use these axioms, except that we allow the functions
from $G^i$ to $V$ to have some sort of singularities.

In practice it is necessary to allow $G$ to be slightly more general than a
group: for example it could be a formal group or Lie algebra.  To
include all these cases together we use cocommutative Hopf algebras $H$.
A typical example of a cocommutative Hopf algebra is the group ring of
a group, and representations of the group are the same as
representations of its group ring. Similarly formal groups and Lie algebras
can also be considered as special cases of Hopf algebras.
In most of this paper $G$ will not be a group, but will be a cocommutative
Hopf algebra together with a ``vertex structure'', as described below.

Informally we can think of any Hopf algebra $H$ as being something
like the group ring of a group, and can think of its dual $H^*$
as being something like the ring of regular functions on this group.
In order to define $G$ vertex algebras we will need to know what is meant
by a ``singular function on $G$'', and to specify these
we need to give a vertex structure on the Hopf algebra $H$.
This means that we are given an algebra $K$ over the ring $H^*$ which
behaves as if it were the ring of singular functions on some group $G$.

The simplest nontrivial example of a vertex structure on a Hopf
algebra is given by taking the Hopf algebra $H$ to be the Hopf algebra
of the one dimensional formal group, and by defining the ring of
singular functions to be the quotient field of the ring of regular
functions. Then the commutative $G$ vertex algebras for this $G$ turn
out to be exactly the same as vertex algebras (as defined in [B] or
[K]). (The (possibly non-commutative) $G$ vertex algebras for this $G$
are also equivalent to previously defined algebraic structures: they
are more or less the same as ``field algebras'' defined in [K] and to
``quantum algebras'' defined in [L-Z].)  In other words vertex algebras
are just commutative $G$ vertex algebras for the simplest possible
nontrivial example $G$ of a Hopf algebra with a vertex structure.

It seems to be hard to construct interesting examples of $G$ vertex algebras 
for higher dimensional $G$, possible because they tend to be closely
related to higher dimensional quantum field theories, which are notoriously
difficult to construct. However we can construct $G$ vertex algebras
corresponding to generalized free field theories, which we do in section 8.

In section 9 we prove the main theorem of this paper, which is a
generalization of the vertex algebra identity [K, 4.8] to $G$ vertex
algebras for higher dimensional $G$. The main result is that (roughly
speaking) the integral of a vertex operator over an $n$-dimensional cycle
is a vertex differential operator of order $n$.

Section 10 and 11 give a few miscellaneous thoughts about $G$ vertex algebras.

There have been several previous definitions of mathematical structures
inspired by quantum field theory; in particular there are
various rigorous definitions of a quantum field theory based
on operators on Hilbert spaces (see [S-W] or [G-J] for example),
and there is Atiyah and Segal's notion of a topological field theory [Se].
In section 12 we discuss the relation between $G$ vertex algebras
and these other concepts.

Beilinson and Drinfeld [B-D] have another approach to defining vertex
algebras, as Lie algebra objects in a suitable category. They remark
that their definition extends to higher dimensions, but the main emphasis 
is on the extension to higher genus curves.  I do not know
what the relation of this is to this paper.  The pseudo tensor
categories in [B-D] are essentially the same as the multilinear
categories of this paper.  Soibelman [So] has some unpublished notes
which  overlap with the present paper.  In particular he
defines ``meromorphic pseudo tensor categories'' which are very
similar to relaxed multilinear categories, and considers associative
algebras in them.  His notes include a generalization to the
``braided'' case, which we largely ignore in this paper.

This is based on a lecture that I promised but did not give at the
1996 Taniguchi symposium. I would like to thank the Taniguchi
foundation for inviting me to Japan, and C. Snydal and
N. R. Scheithauer and the referee for suggesting many corrections and
improvements.

\proclaim 2.~Formal groups and Hopf algebras.

As well as considering groups acting on modules, we also need to consider
various group-like things such as Lie algebras or formal groups
acting on spaces. The most convenient way of doing these is to use
cocommutative Hopf algebras, which are a common generalization of groups, Lie
algebras,
and formal groups.
In this section we recall some standard results about formal groups
and Hopf algebras, most of which can be found in any standard
references such as [A] for Hopf algebras and [S] for formal
groups.

As motivation for the definition of Hopf algebras we consider the
following example.

{\bf Example 2.1.} Suppose $G$ is any discrete group and $H$ is its group ring
over some commutative ring $R$. Then $H$ is an associative algebra
with identity, and also has a map $\Delta$ (called the
comultiplication) from $H$ to $H\otimes_RH$ given by
$\Delta(g)=g\otimes g$ for $g\in G$, and a map $\epsilon$ (called the
counit) from $H$ to $R$ given by $\epsilon(g)=1$, and a map $S$ (called
the antipode) from $H$ to $H$ given by $S(g)=g^{-1}$.  Then $H$ has
the following properties:

\item{1.} $H$ is an associative
algebra with identity over $R$.
\item{2.} $\Delta$ is an algebra homomorphism from $H$ to $H\otimes H$,
which is coassociative, and
$\epsilon$ is a homomorphism  of algebras from $H$ to $R$,
which is a counit for the comultiplication.
\item{3.}
$\mu(S\otimes 1(\Delta(h)))=\mu(1\otimes S(\Delta(h)))=
\epsilon(h)$ for all $h$.

\proclaim Definition 2.2.
A {\bf Hopf algebra} $H$ is defined to be a module over $R$
with multiplication, comultiplication, identity, counit, and antipode $S$,
satisfying the axioms above. If in addition the comultiplication
is cocommutative then the Hopf algebra is called {\bf cocommutative}.

So the group ring of any group is a cocommutative Hopf algebra.
Conversely any cocommutative Hopf algebra behaves in many ways as if it were
the group ring of some group.

The comultiplication is rather difficult to handle in notation.  We
assume that for each $h\in H$ we have chosen finite sets of elements
$h_{1i}$ and $h_{2i}$ such that $\Delta(h)=\sum_ih_{1i}\otimes
h_{2i}$.  If $H$ is a Hopf algebra then the comultiplication on $H$
induces a multiplication on the dual $H^*=\Hom_R(H,R)$ which makes
$H^*$ into an associative algebra over $R$ with unit $\eta$. If $H$ is
cocommutative then $H^*$ is commutative. If $H$ is a finitely
generated free module over $R$ then $H^*$ also has a coalgebra
structure making it into a Hopf algebra induced by the product of $H$,
but in general $H^*$ does not have a coalgebra structure, because the
natural map from $H^*\otimes_RH^*$ to $(H\otimes_RH)^*$ is not usually
an isomorphism.

{\bf Example 2.3.} If $L$ is a Lie algebra over $R$ then the universal
enveloping algebra $U(L)$ of $L$ is a Hopf algebra,
with the comultiplication induced by $\Delta(g)=g\otimes 1+1\otimes g$
and $S$ induced by $S(g)=-g$ for $g\in L\subset U(L)$.
If $R$ is a field of characteristic 0 then Lie algebras
over $R$ are essentially equivalent to formal groups, and the
universal enveloping algebra of $g$ is just isomorphic to
the Hopf algebra of the formal group corresponding to $g$
(defined below). If $R$ is not
an algebra over $\Q$ then the Hopf algebra $U(L)$ usually
does not have particularly good properties, and the Hopf algebras
of formal groups behave much better (at least for the purposes of this paper).

A formal group of dimension $n$ over a commutative ring $R$ (see [S,
II.4.6]) is an $n$-tuple $F(x,y)=(F_1(x,y),\ldots,F_n(x,y))$ of formal
power series in $2n$ variables $x=(x_1,\ldots,x_n),y=(y_1,\ldots,y_n)$
such that

\item{(1)} $F(0,y)=y$, $F(x,0)=x$.
\item{(2)} $F(x,F(y,z))=F(F(x,y),z)$.

Axiom 1 implies the existence of unique power series
$\phi(x)=(\phi_1(x),\ldots,\phi_n(x))$ with $\phi(0)=0$ and
\item{(3)} $F(x,\phi(x))=0=F(\phi(x),x)$.

We will abbreviate $F(x,y)$ to $xy$ and $\phi(x)$ to $x^{-1}$,
so that the axioms above become the same as the usual group axioms
(except that the multiplicative identity is rather confusingly written as 0).

{\bf Example 2.4.} The formal group of $GL_n$. This formal group
has $n^2$ variables $x=(x_{ij})$ which we think of as the $n^2$
entries of an $n\times n$ matrix. The power series $F$ and $\phi$
are defined by $I_n+F(x,y)=(I_n+x)(I_n+y)$ and $I_n+\phi(x)=(I_n+x)^{-1}$
where $I_n$ is the identity matrix. So for example $F_{ik}(x,y)=x_{ik}+y_{ik}
+\sum_jx_{ij}y_{jk}$.

We can also define formal groups of countable infinite dimension in a
countable infinite number of variables $x_i$, $i\in N$, provided we
are careful to use the correct definition of formal power series in an
infinite number of variables: a formal power series in variables $x_i$
is defined as a formal infinite sum $\sum_{0<i_1\le i_2\le\cdots\le
i_n} a_{i_1i_2\ldots i_n}x_{i_1}x_{i_2}\cdots x_{i_n}$. In particular
the ring of formal power series in an infinite number of variables is
strictly larger than the completion of the ring of polynomials in the
$x_i$'s at the ideal generated by the $x_i$'s, because in the latter
there would only be a finite number of nonzero $a$'s for any
given value of $n$.

{\bf Example 2.5.} Suppose $G$ is a formal group of dimension $n$. We can
construct a cocommutative Hopf algebra $RG$, called the formal group
ring of $G$, as follows.  We let $RG$ be a free module over $R$ with a
basis of elements $D^{(i)}$, where $i=(i_1,\ldots,i_n)$ is an
$n$-tuple of nonnegative integers. (If the dimension $n$ is infinite
we should add the condition that all but a finite number of the
$i_j$'s are zero.)  We define the comultiplication by
$\Delta(D^{(i)})=\sum_jD^{(j)}\otimes D^{(i-j)}$, and the antipode $S$
by $S(D^{(i)})=(-1)^{i_1+\cdots i_n}D^{(i)}$. We can identify $RG$ with the
topological dual of the ring of formal power series
$R[[x]]$ by letting $D^{(i)}$ be a dual basis to
the elements $x^i=x_1^{i_1}\cdots x_n^{i_n}$ in the obvious way. Then
the continuous algebra homomorphism from $R[[x]]$ to $R[[x,y]]$ taking
$x$ to $F(x,y)$ induces a dual map from $RG\otimes RG$ to $RG$ which
we use as the algebra structure on $RG$. (Alternatively we could also
define $RG$ to be the algebra of ``continuous left invariant
differential operators on $G$''.) If we work over a field of
characteristic 0, then formal groups are essentially the same as Lie
algebras, and the formal group ring of a formal group is isomorphic to
the universal enveloping algebra of the corresponding Lie algebra.

{\bf Example 2.6.} Suppose $G_a$ is the one dimensional formal group
with formal group law $F(x,y)=x+y$. Then the formal
group ring $H$ over $R$ is the $R$-algebra generated by elements
$D^{(i)}$ for $i$ a non-negative integer,
with $D^{(0)}=1$, $D^{(i)}D^{(j)}={i+j\choose i}D^{(i+j)}$.
If $R$ contains the rational numbers this is just the ring of
polynomials in one variable $D=D^{(1)}$, with $D^{(i)}=D^i/i!$.
The dual $H^*$ can be identified with the ring of Laurent
power series $R[[x]]$, in such a way that
the elements $x^i$ are a dual ``basis'' to the elements $D^{(i)}$.

The category of (left) modules over a Hopf algebra $H$
is defined to be the category of left modules over
$H$ considered as an $R$-algebra. (If $A$ is a right module over $R$
we can turn it into a left module by defining $ha=aS(h)$.)
If $A$ and $B$ are modules over $H$ then
$A\otimes_R B$ and $\Hom_R(A,B)$
can be made into
$H$-modules in a canonical way as follows.
$$h(a\otimes b)=\Delta(h)(a\otimes b)$$
$$(hf)(a) = \sum h_{(1)}(f(S(h_{(2)})(a)))$$
In particular the modules over a Hopf algebra form a tensor abelian
category.

If $P$ is any free module over $R$ then the tensor algebra
$T^*(P)=\oplus_{n\ge 0}\otimes^n(P)$ has a natural coalgebra
structure, given by $\Delta(p_1\otimes\cdots\otimes p_n)=
\sum_{0\le j\le n}(p_1\otimes\cdots\otimes
p_j)\otimes(p_{j+1}\otimes\cdots\otimes p_n)$. We denote the largest
cocommutative subcoalgebra
of this by $B(P)$, so that $B(P)=\oplus_{n\ge 0}(\otimes^nP)^{S_n}$,
where $(\otimes^nP)^{S_n}$ is the subspace of $(\otimes^nP)$ fixed by
the natural action of the symmetric group $S_n$. Over a field of
characteristic 0 the natural map from $B(P)$ to the symmetric algebra
of $P$ is an isomorphism of vector spaces, but this is not true in
positive characteristic. Warning: the symmetric coalgebra $B(P)$ is
not the universal cocommutative coalgebra cogenerated by $P$ (although
some books incorrectly state that it is, probably because of the
analogy with the symmetric algebra). The relation between them is that
$B(P)$ is an irreducible component of the universal cocommutative
coalgebra cogenerated by $P$; see [A].  We can recover $P$ from $B(P)$
as the space of primitive elements, i.e., the elements $p$ with
$\Delta(p)=p\otimes 1+1\otimes p$.
A Hopf algebra is the formal group ring of a formal group over $R$ if
and only if it is a free $R$ module and its underlying coalgebra is
isomorphic to $B(P)$ for some free $R$ module $P$.

\proclaim 3.~Vertex groups.

In this section we define the concept of a vertex group $G$,
which can be thought of informally as a group with a space of functions
with singularities on it.  We
will later define $G$ vertex algebras for any vertex group $G$ as
analogues of algebras acted on by an ordinary group $G$. The definition
of a vertex group has to be weak enough to include
the examples later, and strong enough so that it is
possible to define the concept of an associative algebra over a vertex group.

In the definition below it is helpful to think of $H$ as the Hopf
algebra given by the group ring of some group $G$, so that $H^*$ can
be thought of as something like the ring of functions on the
underlying space of $G$, and the left and right actions of $H$ on
$H^*$ correspond to the left and right actions of $G$ on itself. The
module $K$ below should then be thought of as some sort of ring of
rational or algebraic functions with singularities on $G$.  The
axioms for a vertex structure are based on obvious properties of the
space of singular functions on a group.  In particular the singular
functions form an algebra over the regular functions on which we can
define the operations of left and right translation by elements of the
group and the operation induced by taking inverses of elements of the
group.

\proclaim Definition 3.1. Suppose that
$H$ is a Hopf algebra over a commutative ring $R$.
An {\bf elementary vertex structure} on $H$ consists of an $R$-module $K$ with
the following extra structures.
\item{(1)} $K$ has the structure of an associative algebra over $H^*$.
\item{(2)} $K$ has the structure of
 a 2 sided $H$-module, and the natural map from $H^*$ to
$K$ is a homomorphism of 2 sided $H$-modules. The
product on $K$ is invariant under the left and right actions of $H$.
This means that $h(ab)=\sum (h_{(1)}a)(h_{(2)}b)$ and similarly
for the right action. This axiom can be thought of as saying
that $K$ is closed under ``left and right translation''.
\item{(3)} There is an $R$ linear map $S$ (called the antipode) from $K$ to
$K$ extending the antipode $S$ on $H^*$, and $S(ab)=S(b)S(a)$ whenever each
of $a$ and $b$ is in $H$ or $H^*$ or $K$. This axiom can be thought
of as saying that $K$ is closed under the inversion map of $G$.

\proclaim Definition 3.2.
An {\bf (elementary)  vertex  group} $G$ is a
Hopf algebra $H$
with an (elementary) vertex structure $K$,
such that $H$ is cocommutative, $K$ is commutative, and $S^2=1$.
We call $H$ the group ring of the vertex
group $G$, and we call $K$ the ring of singular functions on $G$.

If we replace the condition $S^2=1$ by the weaker condition that $S$
is an isomorphism of $R$ modules from $K$ to $K$ then we call $G$ a
{\bf braided elementary vertex group}.  There are several other
obvious variations of these definitions (which we will not use in this
paper): if we add the condition that $H$ is the Hopf algebra of a
formal group we get the definitions of vertex formal groups and
braided vertex formal groups, and if we replace the condition that the
vertex structure is cocommutative by the condition that the Hopf
algebra $H$ is ``braided quasi commutative'' or ``quasi triangulated''
we get the concept of a ``braided vertex quantum group''. 

The definition of an elementary  vertex group is not really general enough
for some examples. The main problem is that the only singularities
of functions of several variables we allow are
essentially singularities of functions of two variables,
and there are some examples where we want to allow more
general singularities; see the end of section 8 for examples and
a provisional definition of a non elementary vertex group.
However the definition above is adequate
for most of the examples in this paper (possibly because
most of the examples are related to quantum field theories
that are either free or in small dimensions).
We will usually drop the adjective ``elementary'' from now on.

{\bf Example 3.3.} If $H$ is any cocommutative Hopf algebra, then taking $K=H$
gives $H$ the structure of a vertex group, called the trivial vertex
group structure. In particular vertex groups are generalizations of
cocommutative Hopf algebras (and hence of groups and Lie algebras).

{\bf Example 3.4.} The simplest nontrivial example of a vertex  group
$G$ is given as follows. Suppose that $H$ is the Hopf algebra of the
1-dimensional formal group  $G_a$ of example 2.6, so that $H$ has a basis of
elements $D^{(i)}$ for $i\ge 0$ and $H^*$ can be identified with
$R[[x]]$ as in example 2.6.  We let $K$ be the quotient field
$R[[x]][x^{-1}]$ of $H^*$ consisting of formal Laurent series over $R$, with
$S$
acting as $S(x^i)=(-1)^ix^i$ and $H$ acting as derivations in the
obvious way. (The right and left actions of $H$ on $K$ are identical.)
If $G$ is the vertex group given by $H$ and $K$ then we will
see later that (classical) vertex algebras are exactly the same as the
commutative $G$ vertex algebras as defined in section 6.

{\bf Example 3.5.} Suppose the ring $R$ contains the inverse $1/N$ of some
integer $N$ and has an element $\zeta$ with $\zeta^N=-1$.  Then there
is a variation of example 3.4 where we take $K$ to be the ring
$R[[x^{1/N}]][x^{-1/N}]$ of all formal Laurent series in $x^{1/N}$,
with $S$ acting as $S(x^{1/N}) = \zeta x^{1/N}$.  This gives
examples of braided vertex groups which are not
vertex  groups (because the antipode $S$ does not have period 2).

{\bf Example 3.6.} Take $H$ to be the formal group of the algebraic group
$SL_2$ (over any commutative ring $R$). If we represent an element of
$SL_2(R)$ as ${ab\choose cd}$ then $H^*$ is the ring of formal power
series in $a-1,b,c,d-1$ modulo the ideal $(ad-bc-1)$. We can define a vertex
structure by
putting $K=H^*[b^{-1}]$, and this defines a vertex formal group. Under
the same conditions as example 3.5 we can define a braided vertex formal
group by putting $K=H^*[b^{-1/N}]$.

If $R$ is an integral domain we can often define $K$ to be
the full quotient field of $H^*$, or even the separable algebraic
closure of this quotient field, but this definition of $K$ usually
seems either too large or too small: for most of the examples later in this
paper it is only necessary to invert a few elements of $K$ and inverting
more makes theorem 9.1 much weaker, and on the other hand we also sometimes
want to allow transcendental extensions as in example 3.9 below.

{\bf Example 3.7.} Take $H$ to be the formal group of $SL_n$, so that
$H^*$ is the ring of formal power series in the elements
$a_{ij}-\delta_i^j$ ($1\le i,j\le n$) modulo the ideal generated
by $\det(a_{ij})-1$, where the $a_{ij}$'s are
the entries in an $n$ by $n$ matrix. If $D_k$ is the determinant
of the top right $k\times k$ submatrix, so that $D_n=1$ and $D_1=a_{1n}$,
then we can define a vertex structure by letting
$K$ be $H^*$ with some subset of the $D_k$'s inverted. (So example
3.6 is the case with $n=2$ where we invert $D_1$.)

{\bf Example 3.8.} We can generalize example 3.7 as follows. We let $G$ be a
split simple algebraic group with a fixed choice of Cartan subgroup
$T$ and Borel subgroup $B$. We define $\chi_1,\ldots,\chi_n$ to be the
fundamental characters of $T$ in the Weyl chamber.  These characters
can be extended uniquely to functions on $G$ which are right invariant
under $B$ and left invariant under $\bar B$, which we will again
denote by $\chi_1,\ldots,\chi_n$. For example if $G$ is $SL_n$ and $T$
the diagonal matrices and $B$ the upper triangular matrices then
$\chi_i$ is given by the determinant of the top left $i\times i$
submatrix. Now take $\omega$ to be an element of $G$ representing the
opposition involution of the Weyl group, so that $\omega$ takes $B$ to
$\bar B$, and define functions $D_k$ by $D_k(g)=\chi_k(\omega g)$, so
that these functions are invariant under left and right multiplication
by elements of $B$. We define a vertex structure by inverting
some subset of the elements $D_k$, considered as formal power series
in $H^*$. (Notice that changing $\omega$ only results multiplying
$D_k$ by an invertible power series, so that this vertex
structure does not depend on the choice of representative $\omega$.)
Notice that the product of the functions $D_k$ vanishes exactly
on the complement of the big cell of the Bruhat decomposition,
so that elements of $K$ can be thought of informally
as singular functions defined near the identity of $G$
whose poles are in the complement of the big cell.

{\bf Example 3.9.} Suppose that $R$ is an algebra over the rational numbers and
$H$ is the Hopf algebra
of the formal group of $G_a$ as in example 3.4. Let $K$ be
the ring $R[[x]][x^{-1},\log(x)]$ where $\log(x)$ is a new indeterminate
(written as $\log(x)$ rather than $y$ for mnemonic reasons)
with $D^{(i)}(\log(x))=(-1)^{i-1}x^{-i}/i$ if $i>0$ and $S(\log(x))=
\log(x)$. (We could also define $S(\log(x))=\log(x)+c$ for any $c\in R$
if we did not mind that $S^2\ne 1$.) This defines a vertex
structure on $H$ which gives a vertex group. Similarly in examples 3.6
to 3.8 we could introduce formal logs of the functions $b$ or $D_k$ that
we inverted.

We can define homomorphisms of vertex groups as follows. The definition
is the obvious one if one thinks of a vertex group as a group
ring.
\proclaim Definition 3.10. Suppose $G_1$ and $G_2$ are vertex groups
with group rings $H_1$ and $H_2$ and rings of singular functions
$K_1$ and $K_2$. Then a homomorphism of vertex groups from
$G_1$ to $G_2$ consists of a homomorphism of
Hopf algebras (over $R$) from $H_1$ to $H_2$, together
with a $H_2^*$-algebra homomorphism from $K_2$ to $K_1$
which commutes with the antipode and the left and right actions
of $H_1$.

{\bf Example 3.11.} If $G_1$ and $G_2$ are ``really'' discrete groups, in other
words if $H_1$ and $H_2$ are their group rings
and $K_1=H_1^*$, $K_2=H_2^*$, then a homomorphism of vertex
groups from $G_1$ to $G_2$ is the same as a homomorphism
of groups from $H_1$ to $H_2$.

{\bf Example 3.12.} Suppose $G_1$ is the trivial group (or rather the vertex
group corresponding to it). If $G_2$ is also an honest group, then
there is a unique homomorphism from $G_1$ to $G_2$. If $G_2$ is an
arbitrary vertex group this is not usually true. For example if $G_2$
is the vertex group of example 3.4 with $K_2=R[[x]][x^{-1}]$ then there
is no homomorphism from $G_1$ to $G_2$, because the homomorphism from
$H_2^*=R[[x]]$ to $H_1=R$ cannot be extended to the field $K_2$. The
significance of this is as follows.  If $G_1$ and $G_2$ are groups
then any homomorphism from $G_1$ to $G_2$ induces a forgetful
functor from $G_2$
algebras to $G_1$ algebras in the obvious way. In particular if we
take $G_1$ to be the trivial group this shows (in a rather roundabout
way!) that any $G_2$ algebra has an underlying associative algebra
structure. Similarly for $G$ vertex algebras a homomorphism
from $G_1$ to $G_2$ induces a natural forgetful functor from
$G_2$ vertex algebras to $G_1$ vertex algebras.
However we cannot use this to show that every $G_2$ vertex algebra
has an underlying associative algebra structure, because
the map from the trivial vertex group to $G_2$ need not exist.

{\bf Example 3.13.} We can construct products of vertex groups
with the usual universal properties.  For example
the product $G_1\times G_2$ of two vertex groups $G_1$, $G_2$ is
constructed as follows. The underlying Hopf algebra of $G_1\times G_2$
is $H_1\otimes H_2$, and the ring of singular functions is
$(H_1\otimes H_2)^*\otimes_{H_1^*\otimes H_2^*}K_1\otimes K_2$.

{\bf Example 3.14.} Suppose $G_1$ is the vertex group
with $K_1=R[[z]][z^{-1}]$ of example 3.4 
and $G_2$ is the 2 dimensional vertex group
$G_1\times G_1$, so that $K_2=R[[x,y]][x^{-1},y^{-1}]$.
For each $(a,b)\in R^2$ there is a homomorphism from
$H_1$ to $H_2$ taking $x$ to $az$ and $y$ to $bz$. This extends
to a homomorphism of vertex groups if and only if $a\ne 0$ and $b\ne 0$.

\proclaim  4.~Relaxed multilinear categories.

We can define multilinear maps between the representations of a vertex group,
but the composition of multilinear maps is not in general multilinear.
We deal with this problem by defining relaxed multilinear categories,
where the composition of multilinear maps need not be multilinear,
but can still be compared with multilinear maps.
Symmetric relaxed multilinear categories have the following two properties:
the representations of
a vertex group  form a symmetric relaxed multilinear category,
and it is possible to define algebraic objects like commutative rings
in a symmetric relaxed multilinear category.
The reader willing to assume this can skip the rest of this section,
which consists mainly of content-free category theory.

Soibelman has defined a similar notion in unpublished notes [So],
except that his version is more general in several ways; for example,
he allows braided rather than symmetric categories.

We will start by discussing what is needed in order to be able to
define an associative algebra in some additive category. Obviously it
is sufficient to assume the category is a tensor category, in other
words for every two objects there is a tensor product given with
suitable properties. If we are given good isomorphisms
between $A\otimes B$ and $B\otimes A$ for all $A$ and $B$ then
we get the notions of symmetric and braided multilinear categories.
In any additive symmetric tensor category we can define most common algebraic
structures, such as associative algebras, commutative algebras,
Lie algebras (ignoring problems in characteristic 2), Hopf algebras,
and so on. If $A_1,\ldots, A_n, B$ are objects of any additive tensor category
we have a space of multilinear maps from $A_1,\ldots, A_n$ to $B$,
defined as the maps from $A_1\otimes\cdots\otimes A_n$ to $B$.

We can weaken the definition of tensor category by assuming that we are just
given  spaces of multilinear maps from $A_1,\ldots, A_n$ to $B$ for
objects $A_i$ and $B$, with suitable properties, but are not given 
objects $A_1\otimes\cdots\otimes A_n$ representing them.
(For example, composition is defined under obvious conditions,
and is associative, and the identity maps behave in the obvious way.)
An additive  category with such a structure is called a multilinear category,
and if we add conditions about the action of the symmetric group
on spaces of multilinear maps we get the concepts of symmetric and braided
multilinear categories. In a multilinear category we can still define
associative algebras, commutative algebras, Lie algebras, and so on,
but we cannot define things like coalgebras and Hopf algebras because
these cannot be defined just in terms of multilinear maps but require
maps to tensor products. For a multilinear category to be a tensor category
it is necessary that spaces of multilinear maps should be representable,
but this is not sufficient because there is no reason in general why
the objects $(A\otimes B)\otimes C$, $A\otimes B\otimes C$, and
$A\otimes (B\otimes C)$ should be isomorphic. A multilinear
category with only one object is an operad. 

A multilinear category is sometimes called an additive pseudo tensor category
but this seems a rather misleading name: ``tensor'' suggests that
tensor products exist, and ``pseudo'' then says that in fact they do not.

We now want to define ``relaxed multilinear categories'',
which are a  generalization of multilinear
categories for which the composition of multilinear maps
is not always a multilinear map. 
We start by defining the sieve category $Sieve_n$ as follows.
Its objects are sieves of depth $k\ge 0$ on $1,2,\ldots, n$, by which we mean
a  sequence of equivalence relations $E_0,E_1,\ldots, E_k$
on $1,\ldots, n$ such that
\item{1} $E_k$ is the ``indiscrete'' equivalence relation where any two
elements are equivalent, and $E_0$ is the ``discrete'' equivalence relation
where any element is only equivalent to itself
\item{2} Any equivalence class in any equivalence relation
is an interval $\{i,i+1,\ldots, j-1, j\}$
\item {3} Any equivalence class of $E_{i+1}$ is a union of equivalence
classes of $E_{i}$ (so that the equivalence relations
are increasingly fine).

A sieve can be thought of as a record of someone's attempt
to multiply $n$ non-commuting elements of a ring.
For example suppose someone tries to calculate the product
$a_1a_2a_3a_4a_5a_6$.
At the end of the first day they might have calculated $a_1a_2a_3$ and
$a_5a_6$.
At the end of the second day they might have
worked out $a_1a_2a_3$ and $a_4a_5a_6$, and by the third day they might
have calculated $a_1a_2a_3a_4a_5a_6$. This would correspond to the sieve
$$\eqalign{
E_0&=\{\{1\}\{2\}\{3\}\{4\}\{5\}\{6\}\}\cr
E_1&=\{\{1,2,3\}\{4\}\{5,6\}\}\cr
E_2&=\{\{1,2,3\}\{4,5,6\}\}\cr
E_3&=\{\{1,2,3,4,5,6\}\}\cr
}$$

There are two other useful ways of representing sieves of width $n$
and depth $d$.  The second is as some parentheses, such as
$(((()()()))((())(()())))$, with the property that there are $n$
innermost pairs of parentheses, each of which is contained in exactly
$d$ other pairs of parentheses, and such that there is an outermost
pair containing everything. The pairs of parentheses correspond to
equivalence classes of the equivalence relation, in a way that should
be obvious on comparing the expression above with the sieve in the
example above. For clarity we represent all innermost pairs
of parentheses $()$ by a blob $\bullet$, and miss out the outermost
pair and any other pair which is uniquely determined by
the condition that the depth should be $d$. For example
we would abbreviate the example above to
$(\bullet\bullet\bullet)(\bullet(\bullet\bullet))$.

The third way of representing sieves of width $n$ and depth $d$
is as rooted trees such that all branches are of length $d$,
there are $n$ ``leaves'', and all the edges going upwards
from any given vertex are totally ordered. There is one
vertex of the tree for every pair of parentheses
(or to every equivalence class of every equivalence relation), and an edge
between vertices if one of the vertices corresponds to
a pair of parentheses which is a maximal pair inside the
parentheses corresponding to the other vertex.

The number of sieves of width $n>0$ and depth $d$ is $d^{n-1}$, as it is
easy to see by induction on $n$. 

We say that a sieve $p$ is a refinement of $q$ if every
equivalence relation in $q$ is also in $p$, and we make the sieves
into a category by saying that there is a unique morphism
from any sieve to any refinement. The category of sieves
has an initial element, consisting of the
sieve with just two equivalence relations $E_0$ and $E_1$.

The idea of a relaxed multilinear category is that instead of one
space of multilinear maps from $A_1,\ldots,A_n$ to $B$
we should have several spaces, one for each sieve of size $n$
(so that in some sense the spaces of multilinear maps depend on
``the order in which we calculate the product of elements of the $A_i$'s'').
The space corresponding to a sieve should be related to the
space of any refinement of that sieve.

More precisely,
a relaxed multi category is given by the following data:
\item{(1)} A set of objects.
\item{(2)} For each collection of objects $A_1,\ldots A_n, B$ ($n\ge 1$)
we are given a functor from $Sieve_n$ to sets. The value of
this functor at $p\in Sieve_n$ is denoted by
$Multi_p(A_1,\ldots,A_n;B)$ and is
called the set of multi morphisms from $A_1,\ldots ,A_n$ to $B$ of type $p$.
(It can be thought of as some sort of space of
Taylor series expansions of multilinear maps.)
If $p$ is the initial object of $Sieve_n$ then
we call $Multi_p(A_1,\ldots;B)$ the set of multi morphisms
and sometimes miss out $p$ from the notation.
\item{(3)} If $A_{11},\ldots,A_{n,m_n},  B_1,\dots,B_n,  C$ are objects,
$p_i\in Sieve_{m_i}$ are sieves of the same depth, and   $q\in Sieve_n$, then
we are given a composition map
$$\eqalign{
&Multi_{p_1}(A_{11},\ldots, A_{1m_1}; B_1)\times \cdots\times
Multi_{p_2}(A_{n1},\ldots, A_{nm_n}; B_n)\cr
&\times Multi_q(B_1,\ldots,B_n; C)\cr
\rightarrow &Multi_{q(p_1,\ldots, p_n)}(A_{11},\ldots,A_{nm_n};C)\cr
}$$
where $q(p_1,\ldots,p_n)$ is the sieve given by
first putting the sieves $p_i$ ``side by side'' and then putting the
sieve $q$ ``on top of them''. In terms of parentheses,
this means we replace the $i$'th innermost pair of parentheses
of $q$ by  $p_i$, and in terms of trees we join the trees together by
making the $i$'th leaf of $q$ into the root of $p_i$.

These data should satisfy the following axioms, which we will state
only vaguely as we do not use them later.
\item{(1)} Composition is associative.
\item{(2)} There is an identity morphism in $Multi(A;A)$
for all $A$, with the obvious properties.
\item{(3)} Composition is ``compatible'' with the morphisms
in the categories $Sieve_n$.

If the spaces of multi morphisms  are all abelian groups
and composition is multilinear then we call the category
a relaxed multilinear category.

If we are given isomorphisms from $Multi(A_1,\ldots;B)$
to $Multi(A_{\sigma(1)},\ldots, B)$ for
all $\sigma$ in the symmetric group and these isomorphisms
satisfy various obvious conditions we call the category a relaxed symmetric
multilinear category. Similarly we can define
relaxed braided multilinear categories.

{\bf Example 4.1.} Any additive tensor category or multilinear category
is a relaxed multilinear category, and in these cases
the spaces of multilinear maps do not depend on the choice of sieve $p$.
(In other words the functors from $Sieve_n$ take all morphisms of $Sieve_n$
to the identity morphism.)

We can define associative rings in any relaxed multilinear category as
follows.  An associative ring consists of an object $A$ and maps
$f_n\in Multi(A,A,\ldots A;A) $, $f_1=$ identity, for all $n\ge 1$
(with $n$ copies of $A$ mapping to $A$) such that any composition of
these maps $f_m$ is the image of some $f_{n}$ under the map from
multilinear maps to compositions of multilinear maps.  Informally we
can think of the associativity property in a relaxed multilinear category
like this: the products $(ab)c$ and $a(bc)$ cannot be compared
directly because they lie in different spaces,
but they can both be
compared with $abc$ using the map from the trivial sieve
to the sieves generated by the equivalence relations
$\{\{a,b\}\{c\}\}$ and $\{\{a\}\{b,c\}\}$.

Note that more complicated multilinear maps in a ring in a relaxed
multilinear category are no longer always uniquely determined by
compositions of bilinear maps, so we cannot sensibly define associativity of a
bilinear map in $Multi(A,A;A)$ but can only define associativity of a
sequence of multilinear maps $f_n$ as above.

Similarly we can define commutative rings, Lie algebras, and so
on in any relaxed symmetric multilinear category.
As before, we need to specify all possible products of
bilinear maps (as well as just bilinear maps) as part of the
definitions of these things.

There are several variations of the definition of a relaxed multi
category. For example, instead of using sieves we could use
collections of intervals on $1,2,\ldots, n$ such that and two
intervals in the collection are either disjoint or one contains the
other, and all 1 and $n$ point sets are in the collection. (The union
of all equivalence relations of a sieve is a collection with this
property.) Then composition of multi maps is easier to axiomatize than
if we use sieves, but it is harder to show that the representations of
a vertex group satisfy the axioms.

The difference between an additive tensor category, a multilinear category,
and a relaxed multilinear category can be illustrated
by pretending that all multilinear maps are representable;
we will of course denote the representing objects as tensor products.
In a tensor category the  object $A\otimes B\otimes C$
representing trilinear maps
is isomorphic to $(A\otimes B)\otimes C$.
In a multilinear category these need not be isomorphic
even if both sides exist,
but  we would expect a canonical
map from $A\otimes B\otimes C$ to $(A\otimes B)\otimes C$
corresponding to the composition of two bilinear maps
being a trilinear map. In a relaxed multilinear category
we would expect a map in the other direction from
$(A\otimes B)\otimes C$ to $A\otimes B\otimes C$,
corresponding to the fact that any trilinear map can be related
to  compositions of bilinear maps.

We can illustrate the differences by drawing the diagrams that have
to commute for a bilinear map to be associative. For tensor
categories, the following diagram has to commute:

$$
\matrix{
(A\otimes A)\otimes A &&\cong&&A\otimes (A\otimes A)\cr
\downarrow&&&&\downarrow \cr
A\otimes A&\longrightarrow &A&\longleftarrow &A\otimes A\cr
}$$
For multilinear categories whose multilinear maps are representable
the following diagram has to commute:
$$
\matrix{
(A\otimes A)\otimes A &\longleftarrow&A\otimes A\otimes
A&\longrightarrow&A\otimes (A\otimes A)\cr
\downarrow&&&&\downarrow \cr
A\otimes A&\longrightarrow &A&\longleftarrow &A\otimes A\cr
}$$
For relaxed multilinear categories whose multilinear maps are representable
the following diagram has to commute. Notice that the arrows
in the top row go in the opposite direction from the previous diagram.
Moreover this is only the first of an infinite number of diagrams
describing associativities of more copies of $A$ that we
need to define an associative product.
$$
\matrix{
(A\otimes A)\otimes A &\longrightarrow&A\otimes A\otimes
A&\longleftarrow&A\otimes (A\otimes A)\cr
\downarrow&&\downarrow&&\downarrow \cr
A\otimes A&\longrightarrow &A&\longleftarrow &A\otimes A\cr
}$$

A. Joyal pointed out to me the following similarities between
Stasheff's $A_\infty$ spaces and associative objects in a relaxed
multilinear category.  An $A_\infty$ space is a sort of space with a
product that is associative up to homotopy, and the homotopies are
well defined up to higher homotopies, and so on. More precisely an
$A_\infty$ space $A$ is a space $A$ together with maps $K_n\times
A^n\mapsto A$ for all $n\ge 2$, where $K_n$ is a certain
$n-2$-dimensional cell complex whose cells correspond to topological
types of rooted trees with $n$ leaves together with an ordering on the
branches going upwards from each node. These are almost the same as
sieves as defined above; the only difference is that in a sieve all
the branches have the same height, and this seems to be a minor
technical requirement which could probably be removed by some slight
changes of definition. Moreover the boundary maps in the complexes
$K_n$ correspond to the maps between sieves given by refinement.  By
adjointness we can think of an $A_\infty$ space as a space $A$ with
maps from $A^n$ to $A^{K_n}$ for all $n$ (satisfying various
conditions) which is very similar to the definition of an associative
object in the relaxed multilinear category of representations of a
vertex group $G$, except that we use $A^{K_n}$ instead of the spaces
of singular $A$ valued functions on $G^n$ parameterized by sieves.

\proclaim 5.~The representations of a vertex group.

In this section we will describe how to make the representations of a
vertex group $G$ into a symmetric relaxed multilinear category, so
that it is possible to define commutative ring objects in this
category.  We let $G$ be a vertex group with group ring $H$ and ring
of singular functions $K$.

The category $C$ of representations of $G$ is defined to be
the category of representations of the group ring $H$. (Of course
the multilinear maps in these two categories will usually
be different!) In particular
the category $C$ has tensor products and Hom functors,
as it is just the representations of some cocommutative Hopf algebra. There
is a functor $\Gamma$ from $C$ to $R$-modules taking any
$G$ module to its fixed point set. Note that
$\Hom_R(x,y)$ is an object of $C$, while $\Hom_G(x,y)$
is the morphisms from $x$ to $y$ which is just an $R$ module.
The relation between them is $\Hom_G(x,y)=\Gamma(\Hom_R(x,y))$.

We will first define the multilinear maps
from $A_1,\ldots, A_n$ to $B$ in $C$. To motivate this
we first describe the $H$-multilinear maps in a slightly unusual way.
The $R$-multilinear maps are just the elements of
$\Hom_R(A_1\otimes\cdots\otimes A_n,B)$,
and the $H$ multilinear ones are just the $H$ invariant elements of
this space. Now consider the map
taking $a_1\in A_1,\ldots, a_n\in A_n$
to $a_1^{g_1}a_2^{g_2}\cdots a_n^{g_n}\in B$ for
$g_i\in H$. For fixed $a_i$ this can be thought of
as a $G$ invariant $B$ valued function on $G^n$, by which we mean an
element of $\Hom_H(H\otimes \cdots\otimes H,B)$. Moreover
this map from $A_1\otimes\cdots\otimes A_n$ to $B$ valued functions
on $G^n$ is $G^n$-invariant. So to summarize, the multilinear
maps can be thought of as $G^n$ invariant maps from
$A_1\otimes\cdots\otimes A_n$ to the $G$-invariant $B$-valued functions
on $G^n$.

We will define $G$ multilinear maps in the same way, except
that we allow the functions on $G^n$ to have ``singularities''.
Informally we want this space to be the space of functions on
$G^n$ which are allowed to have singularities ``of type $K$'' whenever
two of the components of $G^n$ are equal. The formal definition
goes as follows.
The space $\Hom(H^n,B)$ of functions from $G^n$ to $B$
is a module over the ring $H^{n*}$ of functions on $G^n$. For each
$1\le i<j\le n$ there is a homomorphism from $H^*$ to $H^{n*}$
induced by the map from $H^n$ to $H$ taking $g_1\otimes\cdots \otimes g_n$
to $g_iS(g_j)$ (= ``$g_ig_j^{-1}$''). We define  the localization
of a module $M$ over $H^{n*}$ at $(i,j)$ to be
$M\otimes_{H^*}K$, where $M$ is made into an $H^*$ module
using the homomorphism $f_{ij}$ from $H^*$ to $H^{n*}$.
\proclaim Definition 5.1.
We define the space $Fun(G^n,B)$ of singular functions from $G^n$ to $B$
to be the localization at all $(i,j), 1\le i<j\le n$ of
the space of functions from $G^n$ to $B$.

The module $\Hom(H^n,B)$ also appears in the computation
of the cohomology of the $G$-module $B$ using the homogeneous
standard resolution: the cohomology groups of $B$
are just the cohomology groups of a complex whose
terms are the $G$-invariant elements $\Hom(H^n,B)^G$. We will
not use this connection in this paper.

{\bf Example 5.2.} Suppose $G$ is the vertex formal group of example 3.4,
so that $H^* =R[[x]]$ and $K=R[[x]][x^{-1}]$. Then the
functions from $G^3$ to $B$ can be identified with the
space of formal power series $B[[x_1,x_2,x_3]]$ in 3 variables
with values in $B$. The localization of this at $(i,j)$ for
$1\le i<j\le 3$ is just $B[[x_1,x_2,x_3]][(x_i-x_j)^{-1}]$,
so the space of singular functions from $G^3$ to $B$ is the space
$$B[[x_1,x_2,x_3]][(x_1-x_2)^{-1},(x_2-x_3)^{-1},(x_1-x_3)^{-1}].$$

\proclaim Definition 5.3. The $G$-module $Multi^K_R(A_1,\ldots,A_n,B)$
of multilinear maps
is defined to be
the space of $H^{n*}$ invariant maps from
$A_1\otimes\cdots\otimes A_n$ to the space of singular $B$-valued
functions on $G^n$.
The $R$-module $Multi^K_G(A_1,\ldots,A_n,B)$ of $G$ invariant
singular multilinear maps  is defined
to be the module of $G$-invariant elements of $Multi^K_R(A_1,\ldots,A_n,B)$.

We have defined the 4 different spaces of multilinear maps
$Multi_R(A_1,\ldots,A_n,B)$, $Multi_G(A_1,\ldots,A_n,B)$,
$Multi_R^K(A_1,\ldots,A_n,B)$, and $Multi_G^K(A_1,\ldots,A_n,B)$.
The first two are the ``usual'' spaces of multilinear and
invariant multilinear maps over the Hopf algebra $H$,
and the other two are the corresponding spaces
of singular multilinear maps.

To make $C$ into a relaxed multilinear category we not only have to
define the multilinear maps, but also the spaces $Multi_p(A_1,\ldots,
B)$ for nontrivial sieves $p$, and the maps between them for $p<q$,
and the compositions of multilinear maps.

We first define the multilinear maps of type $p$, for $p$
in the sieve category $Sieve_n$. These can be thought of informally
as various sorts of Taylor series expansions of multilinear maps.

\proclaim Definition 5.4. The
space $Fun_p(G^n,B)$ of singular $B$-valued functions  of type
$p$ is defined recursively as follows. Represent $p$ as a tree
of depth $d$, consisting of $n$ leaves attached to
a tree $q$ of depth $d-1$. Suppose that $X$ is the
space of singular $B$-valued functions of type $q$.
Then $Fun_p(G^n,B)$ is obtained from $X$ in two steps: first we construct
the space $Y$ of all $G^n$ valued functions with values in $X$,
then we localize $Y$ at all pairs $i,j$ for $i$ and $j$ leaves
of $p$ joined to the same leaf of $q$.

This space can be thought of as some space of singular functions 
``of type $p$'' on $G^n$
with values in $B$, where $n$ is the number of nodes of the tree $p$.
This space is also acted on by $G^n$.

{\bf Example 5.5.} With $G$ as in example 3.4, we have
$$\eqalign{
Fun_{(\bullet\bullet\bullet)}(G^3, B)
&=B[[x_1,x_2,x_3]][(x_1-x_2)^{-1},(x_2-x_3)^{-1},(x_1-x_3)^{-1}]\cr
Fun_{((\bullet\bullet)\bullet)}(G^3, B)
&=B
[[x_1,x_2]][(x_1-x_2)^{-1}]
[[x_3]][x_3^{-1}]\cr
Fun_{(\bullet(\bullet\bullet))}(G^3, B)
&=B
[[x_2,x_3]][(x_2-x_3)^{-1}]
[[x_1]][x_1^{-1}]\cr
}$$

\proclaim Definition 5.6. The $G$-module $Multi^K_{p}(A_1,\ldots, A_n, B)$
of singular multilinear maps of type $p$ is defined
to be
$$\Hom_{H^{n}}(A_1\otimes\cdots\otimes A_n, Fun_p(G^n,B)).$$
The $R$-module $Multi^K_{Gp}(A_1,\ldots, A_n, B)$ of $G$-invariant
singular multilinear maps of type $p$ is defined to be the
$G$-invariant elements of $Multi^K_{p}(A_1,\ldots, A_n, B)$.

The composition of multilinear maps of types $p_1,\ldots, p_n$ and $q$
is easily checked to be a multilinear map of type $q(p_1,\ldots, p_n)$.

Finally we have to describe the map from multilinear maps of type
$p$ to multilinear maps of type $q$ whenever $p<q$.
The only possible problem in defining this map is showing
that the localizations  extend. The way to do this is best explained by an
example.
If we have an expression like
$(x_1-x_2)^{-k}$
we can expand it as a series
$$\sum_{i,j\ge 0}C(i,j) (x_1-y_1)^{i}(y_1-y_2)^{-k-i-j}(y_2-x_2)^{j}$$
for some constants $C(i,j)=(-k)!/i!j!(-k-i-j)!$.
In general there is an analogue of this power series expansion
for any cocommutative Hopf algebra $H$.

\proclaim 6.~$G$ vertex algebras.

In this section we give the definition of a $G$ vertex algebra.

We start with some motivation for the definition. Suppose that $V$ is
an associative algebra acted on by a group $G$, with the action of
$g\in G$ on $v\in V$ written as $v^g$.  For each fixed $v_1,\ldots,v_n\in V$
we can form the products $v_1^{g_1}v_2^{g_2}\cdots v_n^{g_n}$, which
we think of as a function from $G^n$ to $V$. These function have the following
properties:
\item{1.} Identity: $\cdots v_n^{g_n}1^gv_{n+1}^{g_{n+1}}\cdots = \cdots
v_n^{g_n}v_{n+1}^{g_{n+1}}\cdots$
\item{2.} Associativity and $G$ invariance:
$$(v_{11}^{g_{11}}v_{12}^{g_{12}}\cdots)^{g_1}
  (v_{21}^{g_{21}}v_{22}^{g_{22}}\cdots)^{g_2}
\cdots
= v_{11}^{g_{11}g_1}v_{12}^{g_{12}g_1}\cdots
  v_{21}^{g_{21}g_2}v_{22}^{g_{22}g_2}\cdots
$$

Conversely if we have a module $V$ over a ring $R$ together
with functions from $G^n$ to $V$ as above which are multilinear
in the $v_i$'s then we can define the structure of an associative
$G$ algebra on $V$ by defining the product $v_1v_2$ to be $v_1^{1}v_2^1$
and defining the $G$ action on $V$ by
$v^g$. It is easy to check that this means that modules $V$ with
a set of functions as above are equivalent to associative algebras
with an action of $G$.

We now define a $G$ vertex algebra where $G$ is a vertex group.  We
defined a vertex group earlier as a cocommutative Hopf algebra over
$R$ with a vertex structure, but we can think of it informally as
being a group together with some ring of ``singular'' functions on it.
A $G$ vertex algebra is informally defined to be a module over $R$
together with multilinear singular maps denoted by
$v_1^{g_1}v_2^{g_2}\cdots v_n^{g_n}$ from $G^n$ to $V$ for all
$v_1,\ldots,v_n\in V$, such that the axioms above are satisfied. In
other words we just copy the
strange definition of an associative $G$-ring in the paragraph above,
and replace functions on $G$ with singular functions on $G$.  The
formal definition, which makes precise what we mean by a singular
function, is as follows.
\proclaim Definition 6.1. A $G$ vertex algebra is
an associative algebra in the relaxed multilinear category of
$G$-modules.

If we give $G$ the trivial vertex structure such that the singular
$V$-valued functions on $G^n$ are defined to be the ordinary functions
$\Hom(G^n,V)$ on $G^n$ then we see from the previous remarks that a
$G$ vertex algebra is exactly the same as an associative algebra acted
on by $G$ (although defined in a rather odd way).

It is now easy to extend many concepts in ring theory
to $G$ vertex algebras as in the following examples. 

{\bf Example 6.2.} A $G$ invariant left ideal of a $G$ vertex
algebra $V$ is an $R$ submodule $I$ such that $v_1^{g_1}v_2^{g_2}\ldots
v_n^{g_n}$
is a singular map from $G^n$ to $I$, whenever $v_1,\ldots,v_{n-1}\in V$,
$v_n\in I$.
Similarly we can define right and two sided ideals.

{\bf Example 6.3.} A left module over the $G$ vertex algebra $V$ is an $R$
module $M$ together with multilinear singular maps $v_1^{g_1}\cdots
v_n^{g_n}m$ from $G^n$ to $M$ for all $v_i\in V, m\in M$, satisfying
the obvious axioms about associativity and the identity element. This
is the definition of a module without an action of $G$; we can define
modules with an invariant action of $G$ in the same way by using
multilinear singular maps $v_1^{g_1}\cdots v_n^{g_n}m^g$ from
$G^{n+1}$ to $M$ instead. We can define submodules, homomorphisms,
and multilinear maps of modules in the obvious way. Unlike the case
of rings there seems to be no reason in general why multilinear maps
should be representable (though they often are in practice),
but when they are representable we can use this to define tensor
products of modules. (The problem of when this tensor product
is associative seems rather hard.)

{\bf Example 6.4.} We can define a commutative $G$ vertex algebra
to be one such that the multilinear singular functions
are invariant under the obvious action of the symmetric group.
We can define supercommutative $G$ vertex algebras in a similar way.

{\bf Example 6.5.} If $V$ and $W$ are two $G$ vertex  algebras
then their tensor product is also a $G$ vertex algebra,
with the $G$ vertex algebra structure defined by
$$(v_1\otimes w_1)^{g_1}\cdots
= (v_1^{g_1}\cdots)\otimes (w_1^{g_1}\cdots).$$
If $V$ and $W$ are commutative then so is $V\otimes W$.
(Note that if $V$ and $W$ are supercommutative
then there should be some extra signs $(-1)^{\deg(v_i)\deg(w_j)}$
for $i<j$ added to the definition of the product in $V\otimes W$.)
More generally, if $V$ and $W$ are $G_1$ and $G_2$ vertex algebras, 
then $V\otimes W$ is a $G_1\otimes G_2$ vertex algebra. 
If $G_1=G_2=G$ then there is a homomorphism from $G$ to $G\otimes G$,
which can be used to turn $V\otimes W$ into a $G$ vertex algebra
as in the first part of this example. 

{\bf Example 6.6.} If $G$ is the vertex group of example 3.4
then it is easy to check that commutative $G$ vertex algebras are 
the same as vertex algebras (as defined in [K, 1.3] for example). 
In fact if $V$ is a vertex algebra and $v_1,\ldots,v_n$ are any elements of
it, and $v_i(z)$ is the vertex operator of $v_i$, then 
$v_1(z_1)\cdots v_n(z_n)1$ is an element of 
$$V[[z_1,\ldots,z_n]][\prod_{i<j}(z_i-z_j)^{-1}]$$ and hence defines a
singular multilinear map from $V\times\cdots\times V$ to $V$, and it
is not hard to check that the axioms for a vertex algebra in [K, 1.3]
are equivalent to saying that this sequence of multilinear maps makes
$V$ into a commutative $G$ vertex algebra.  The field algebras of [K,
4.11], which are ``non-commutative'' vertex algebras, seem to be the
same as associative $G$ vertex algebras for $G$ as in example 3.4.

{\bf Example 6.7.} If $G$ is the vertex group whose underlying Hopf
algebra is the universal enveloping algebra of the Virasoro algebra, with the
ring of singular functions given by inverting an element as in example
3.4, then $G$ vertex algebras are more or less the same as vertex
operator algebras, provided that we restrict to the positive energy
representations of the Virasoro algebra.

{\bf Example 6.8.} If $G$ is the product of two copies of the
vertex group of example 3.4, then $G$ can be thought of
as a vertex group corresponding to 2 dimensional quantum field
theories. A typical $G$ vertex algebra can be constructed by
taking the tensor product of 2 vertex algebras; this corresponds
to the usual decomposition of fields in 2 dimensional field theory
into their left moving and right moving parts. We can do the same thing
if $G$ is the product of 2 copies of the vertex group in example 6.7;
this gives the vertex group underlying much of conformal field theory 
and string theory, and as before we can construct many examples 
by taking a tensor product of 2 vertex operator algebras.

One reason why nontrivial examples of $G$ vertex algebras
are  hard to find is that
there are almost no examples of finite dimensional $G$ vertex algebras,
except for those that are really just associative algebras.  It is not
hard to see why this should be so: a nontrivial $G$ vertex algebra
should have some sort of singularity, say a pole of order 1, otherwise
it would just be an associative algebra.  However from a pole of order
1 we can usually construct poles of all orders by algebraic
operations, so we get an infinite dimensional space of possible
singularities. This is not possible unless the $G$ vertex algebra we
started with is infinite dimensional.

\proclaim 7.~Vertex operators.

The most powerful and general way of constructing $G$ vertex algebras
is to write down a set of vertex operators satisfying the ``normal
ordering'' condition below. This is closely related
to the well known method of constructing vertex algebras by writing
down a set of commuting vertex operators; see [K, 4.5]. 
For the rest of this section we suppose that $G$ is a vertex  group
with group ring $H$. We now define several rings of
``singular functions'' and some modules over them.

If $V$ is an $R$ module we define $V[[x]]$ to be $Hom_R(H,V)$ (where
$x$ stands for $\dim(G)$ variables). If $G$ is a formal
group this is isomorphic as an $R$
module to the space of power series in $\dim(G)$ variables with
coefficients in $V$, and is an $H$ module (using the trivial action of
$H$ on $V$). Similarly we define $V[[x_1,\ldots,x_n]]$ to be
$Hom_R(H\otimes\cdots\otimes H,V)$, which is isomorphic to the space
of power series in $n\dim(G)$ variables with coefficients in $V$ and
is an $H\otimes\cdots\otimes H$ module.  It is also a module over the
ring $R[[x_1,\ldots x_n]]=
\Hom(H\otimes\cdots\otimes H,R)$.
We define $V((x_1,\ldots,x_n))$ to be $V[[x_1,\ldots,x_n]]$ localized
at all $1\le i\le n$ and at all pairs $1\le i<j\le n$ as in section 5.
For example, if $G$ is as in example 3.4 then $$V((x_1,x_2))=
V[[x_1,x_2]][x_1^{-1}x_2^{-1}(x_1-x_2)^{-1}].$$ 

\proclaim Definition 7.1. A vertex operator on an $R$ module $A$ is an
$R$ linear map from $ A $ to $B((x))=Hom_R(H,B)\otimes_{H^*}K$.
More generally a vertex
operator from $A_1\times\cdots\times A_n $ to $B$
is an $R$ multilinear map from $A_1\times\cdots\times A_n $ to
$B((x_1,\ldots,x_n))$.

{\bf Example 7.2.} With $G$ as in example 3.4, a vertex operator
from $A$ to $B$ is just a linear map from $A$ to the space of formal
Laurent series in $B$.

A vertex operator on $A$ can be thought of informally as
a singular function from $G$ to operators from $A$ to $B$.

It is not usually true that the composition of vertex operators
is a vertex operator.
\proclaim Definition 7.3.
We say that a sequence of vertex operators $v_1,\ldots,v_n$
is {\bf compatible} if their composition is
the image of a unique singular multilinear map.
More precisely suppose each $v_i$ is a vertex operator from
$A_{i-1}$ to $A_i$, in other words an $R$ linear map
from $A_{i-1}$ to $A_i((x_i))$. Then the composition is
automatically a linear map from $A_0$ to $A_n((x_n))\cdots((x_1))$,
and the condition that the composition is the image
of a unique singular multilinear map  means
that it is in the image of a unique map from
$A_0$ to $A_n((x_n,\ldots,x_1))$.
We say that a  set $S$ of vertex operators on $A$
is compatible if any finite sequence of vertex operators in
$S$ is compatible.

\proclaim Definition 7.4.
We say that a vertex operator from $A$ to $B$ is a {\bf creation
vertex operator} if it is (induced by) a linear map from $A$ to
$B[[x]]=Hom(H,B)$.

\proclaim Definition 7.5. We say that a vertex
operator from $A$ to $B$ is an {\bf annihilation vertex operator} if
it factors through $B\otimes_RK$ via the natural map from
$B\otimes_RK=B\otimes_RH^*\otimes_{H^*}K$ to $\Hom(H,B)\otimes_{H^*}K$.

Roughly speaking, a creation vertex operator is
a vertex operator that is nonsingular, and
an annihilation vertex operator is one that is
``rational'' rather than ``transcendental''.

\proclaim Definition 7.6. We say that a product of
vertex operators is {\bf normally ordered}
if for some $k\ge 0$ the first $k$ operators in the product
are creation vertex operators and the remainder are
annihilation vertex operators.

\proclaim Definition 7.7.
We say that a set of vertex operators $S$ on an $R$-module $A$
satisfies the {\bf normal ordering condition} if the composition of
any two vertex operators in $S$ can be written as a linear combination
of normally ordered products of pairs of vertex operators in $S$.

\proclaim Lemma 7.8. If a set of vertex operators
satisfies the normal ordering condition
then it is a compatible set of vertex operators.

Proof. The main point is that any normally ordered product of
annihilation and creation operators is a vertex operator.
It follows immediately that if a set of vertex operators satisfies
the normal ordering condition then any composition of vertex operators in $S$
can be written as a linear combination of normally ordered
compositions of vertex operators in $S$, by repeatedly rewriting product
of operators so that creation vertex operators occur on the left, and
hence is a vertex operator.  So the set $S$
is a compatible set of vertex operators. This proves lemma 7.8.

The  point of the preceding definitions is the following theorem,
which will be the main tool for finding compatible sets of vertex operators.
It is an analogue of the trivial theorem that any set of operators
on a module generates an algebra acting on the module.
\proclaim Theorem 7.9. Any compatible set of vertex operators
on a module $A$ generates a $G$ vertex algebra acting on $A$.

If $A$ is a module over a commutative ring $V$ which is generated
as a $V$ module by an element $1_A$, then $A$ obviously has
the structure of a commutative ring, isomorphic to the kernel
of $V$ by the ideal of elements annihilating $1_A$. The next theorem
is a generalization of this to vertex algebras.

\proclaim Theorem 7.10. Suppose an $R$  module $A$ is acted on by a commutative
$G$ vertex algebra $V$, and suppose that $A$ contains an element $1_A$
fixed by $G$ and such that $A$ is generated as a $V$-module by $1_A$
(which means that there is no proper $G$-submodule of $A$ containing $1_A$
and acted on by $V$) and such that
the vertex operator of any element of $V$ applied to $1_A$ is
nonsingular. Then the module $A$ has a unique structure
of a commutative vertex algebra such that the identity is $1_A$
and such that the action of any element of $V$ on $A$
is given by the action of some element of $A$.

Proof. We define a map from $V$ to $A$ by mapping
any $v$ to $v^1(1_A)$ (which is well defined as we assumed
$v^g(1_A)$ to be nonsingular in $g$). This is a homomorphism
of $V$ modules, and as $V$ is commutative the kernel is
a two sided ideal of $V$, so the image of $V$ is a commutative
vertex algebra. As this image is a $V$ submodule containing
$1_A$, it must be the whole of $A$ by assumption. Therefore
$A$ has the structure of a commutative $G$ vertex algebra.
This proves theorem 7.10.

The usual method of using the theorems in this section to construct
$G$ vertex algebras goes as follows. First we write down
a lot of annihilation and creation operators
satisfying the normal ordering condition. Then we try to find
a commuting set of vertex operators in the vertex algebra they generate. There
is often then an obvious element $1_A$ satisfying the conditions
of theorem 7.10, so this makes $A$ into a commutative $G$ vertex algebra.

{\bf Example 7.11.} We show how to construct the vertex algebra of a lattice
using theorem 7.10. (See [K, 5.4] for a similar construction.)  We
assume for simplicity that all inner products in the lattice $L$ are
even. We let $G$ be the vertex group of example 3.4 with underlying
Hopf algebra $H$ the ring spanned by the elements $D^{(i)}$. We let
$V$ be the universal commutative $H$-algebra generated by the group
ring $R[L]$ of the lattice $L$. For each $\alpha$ in $L$ we define a
creation operator $e^{\alpha+}(z)$ to be multiplication by the element
$\sum z^iD^{(i)}(e^\alpha)$.  We define the annihilation operator
$e^{\alpha-}$ to be the homomorphism of rings with derivation from $V$
to $V[z,z^{-1}]$ taking $e^\beta$ to
$z^{(\alpha,\beta)}e^\beta$. These annihilation and creation operators
generate a (non commutative) $G$ vertex algebra acting on $V$, by
theorem 7.9. The operators
$e^{\alpha}(z)=e^{\alpha+}(z)e^{\alpha-}(z)$ all commute with each
other, and applying theorem 7.10 shows that they induce the structure
of a vertex algebra on $V$. This is the usual vertex algebra
of a lattice $L$. 

\proclaim 8.~Free field $G$ vertex algebras.

In this section we construct some examples of $G$ vertex algebras
closely related to generalized free quantum fields.

To keep notation simple we will first describe a special case of the
construction with just one (scalar) field. Afterwards we will list
various ways to generalize it.  We assume that the underlying Hopf
algebra of $G$ is the universal enveloping algebra of the Lie algebra
of the group of
translations of spacetime, so it is a polynomial algebra generated by
the elements ${\partial\over \partial x_i}$ for $1\le i\le d$. The dual
$H^*$ is the ring of formal power series $R[[x_1,\ldots,x_d]]$.
We will take $K$ to be  $H^*[(-x_1^2+\cdots+x_n^2)^{-1}]$ 
(so we allow functions to have singularities along the light cone).  We
will say that an element $f$ of $H^*$ or $K$ is even if $S(f)=f$ and
odd if $S(f)=-f$, where $S$ is the antipode with $S(x_i)=-x_i$.
We fix an even element
$\Delta(x)$ of $K$, which we call the propagator.

We let $A$ be a free module over $H$ generated by
an element $\phi$ (which will be the free field).
The underlying space of the $G$ vertex algebra $V$ we are constructing
is the symmetric algebra $S^*(A)$. This is an associative commutative
algebra acted on by $G$, and we will call the
product of elements $a$ and $b$ in this algebra the normal ordered product
of $a$ and $b$ and denote it by $:ab:$

We will construct an annihilation
vertex operator $\phi^-(x) $ and a creation vertex operator
$\phi^+(x)$, and we  define $\phi(x) $ to be $\phi^-(x)+\phi^+(x)$.
We define the annihilation operator $\phi^-(x)$  to
be the derivation of $R[[x]]$ algebras from $V[[x]]$ to $V[[x]]$ which
commutes with the action of $G$ and such that
$$\phi^-(x)(\phi) = \Delta(x).$$
We define the creation operator from $V$ to $V[[x]]$ by
$$\phi^+(x)(v)= \sum_ix^iD^{(i)}(\phi)v.$$

\proclaim Lemma 8.1. The annihilation and creation operators satisfy
the following commutation relations.
$$\eqalign{
[\phi^+(x),\phi^+(y)]&=0\cr
[\phi^-(x),\phi^+(y)]&=\Delta(x-y)\cr
[\phi^-(x),\phi^-(y)]&=0\cr
[\phi (x),\phi(y)]&=0\cr
}$$

Proof. The first equality is trivial because the ring $V$ is commutative,
and the fourth follows from the first three and the fact that
$\Delta$ is even. The third equality
follows because $[\phi^-(x),\phi^-(y)]$ is a commutator of 2 derivations
on $V[[x,y]]$ and is therefore also a derivation, and it
commutes with $G$ and vanishes on $V$
 so it is zero.
For the second equality we calculate that
$$\eqalign{
\phi^-(x)(\sum_jy^jD^{(j)}(\phi))
&=\sum_jy^jD^{(j)}\Delta(x)\cr
&=\Delta(x-y)\cr
}$$
which implies the second equality  because $\phi^-(x)$ is a derivation.
This proves lemma 8.1.

\proclaim Theorem 8.2. There is a unique structure of a commutative
$G$ vertex algebra
on $V$ such that the vertex operator $\phi(x)$
is the vertex operator of some element.

Proof. This follows easily from lemma 8.1 and theorem 7.10 and some
routine checks.

The elements of the $G$ vertex algebra $V$ are just the usual
Wick polynomials in $\phi$, and the vertex algebra
product of elements of $V$ is just the usual expansion of
products of Wick polynomials in terms of normal ordered products.
See [S-W], section 3.2, for examples.

The construction above can easily be generalized in several ways as
follows.
\item{1.} We can allow more than one field $\phi$,
in which case $\Delta$ should be changed to
a function $\Delta_{\phi\psi}$ depending on the fields
$\phi$ and $\psi$.
\item{2.} We can enlarge $G$ to a semidirect product
of spacetime translations with some other group,
such as the Lorentz group. This other group will act on the space
of fields in 1 above, so we can define spinor fields, vector fields,
and so on.
\item{3.}If $\Delta$ satisfies some differential equation
(such as the wave equation),
in other words if it is annihilated by some element
$D$ of $H$, then $D(\phi)$  generates a proper ideal of the
$G$ vertex algebra $V$ so we can quotient out by it.
\item{4.}Suppose that we have a function
$\Delta_n(g_1,g_2,\ldots,g_n)$ of $n$ variables $g_i\in G$
which are $g$ invariant, with $\Delta_n$ nonsingular for $n>2$.
Then we can define a new $G$ vertex algebra structure on
$S^*(A)$ so
that $\phi(x_1)\phi(x_2)\cdots \phi(x_n)$
to be $:\phi(x_1)\phi(x_2)\cdots \phi(x_n):+\Delta(x_1,\cdots, x_n)$
More generally we define products of the form
$\phi^{n_1}(x_1)\phi^{n_2}(x_2)\cdots $
to be a sum of terms, each of which is formed by repeatedly pulling
out $n$ factors of the form $\phi(x_i)$ with not all $i$'s the same,
and replacing them with $\Delta_n$.
Derivatives of $\phi$
are handled by differentiating $\Delta$.
More generally we can assume that we are given functions $\Delta_n$
for all $n\ge 2$ as above.
The functions $\Delta_n$ are the ``irreducible $n$-point functions''.
\item{5.} In example 4
we need to assume that the functions $\Delta_n$ are nonsingular
for $n>2$ in order to fit into the framework of $G$ vertex algebras.
We can allow the functions $\Delta_n$ to be singular
(but we have to assume that if we identify some but not all of
the variables of $\Delta$ with each other then the
result is a well defined function) if we are willing
to enlarge the notion of a vertex group.
In particular we have to allow singular multilinear functions
to have more general sorts of singularities, so the more general definition
of a vertex group should specify not just the singular functions
of one variable, but also singular functions of many variables.
We will leave the precise definition as an exercise for the reader.

Example 5 above gives examples corresponding to nontrivial
quantum field theories defined perturbatively, which corresponds to the
fact that the $G$ vertex algebras are only defined over
the ring of formal power series in the coupling constants.

\proclaim 9.~The main identity.

In this section we prove an identity for commutative 
$G$ vertex algebras where $G$ is a vertex formal  group associated to
some algebraic groups over $\C$.  Roughly speaking, the identity says
that if we have a vertex differential operator acting on a commutative
vertex algebra, then integrating it over a cycle of dimension $n$
increases its order as a differential operator by at most $n$.

In the case of classical vertex algebras (when $G$ as in example 3.4
is associated to the one dimensional additive algebraic group) this
identity we will prove is equivalent to the identities originally used
to define vertex algebras [B] (see [K]). We will start off by
explaining why, as motivation for the proof.
 The vertex algebra identity states that
$$\sum_{i\in \Z}{m\choose i}(u_{q+i}v)_{m+n-i}w
 = \sum_{i\in \Z}(-1)^i{q\choose i}(u_{m+q-i}(v_{n+i}w)
-(-1)^qv_{n+q-i}(u_{m+i}w)).$$
For simplicity we will only discuss the
case $m=n=q=0$, when it becomes (one version of) the
Jacobi identity $(u_0v)_0w= u_0(v_0w)-v_0(u_0w)$;
the general case follows from a similar argument after including
a suitable rational function of $x$, $y$ and $z$.
This Jacobi identity can be deduced from the fact that $u_0$
is a differential operator of degree 1 as follows. The fact that
$u_0$ is a differential operator of degree 1 means that
$(u_0)v(y)-v(y)u_0$ is a differential operator of degree
0 and must therefore be of the form $t(y)$ for some $t\in V$; putting
$y=0$ then shows that $t=u_0v$, so that $u_0(v(y)w)-v(y)u_0w=(u_0v)(y)w$,
and now integrating $y$ around $0$ gives the Jacobi identity above.
This justifies the statement that the main identity for
a vertex algebra is equivalent to the fact that integrating
$u(x)$ along a 1-cycle is a vertex differential operator of degree 1.

To prove that the integral $u_0$ of $u(x)$ along a 1-cycle
is a vertex differential operator of degree 1
we have to show that the double commutator
$$[[u_0,v(y)]w(z)]$$
is 0, and each term in the double commutator can be written
as the integral of $u(x)v(y)w(z)$ along a suitable cycle in the
$x$ plane with the points $0$, $y$, and $z$ removed. More precisely,
if $C_{a,b,c\ldots}$ is a 1-cycle going once clockwise around
the points $a,b,c\ldots$ and not containing other elements of the points
$0,x,y$ then
$$C_{0,y,z}-C_{0,y}-C_{0,x}+C_0=0$$
in the first homology group, which gives us the relation
$$(u_0)v(y)w(z)-w(z)(u_0)v(y) -v(y)(u_0)w(z)+v(y)w(z)u_0=0.$$
In other words, the fact that $u_0$ is a vertex differential operator
follows from linear relations between elements of a certain homology group.

In the higher dimensional case we proceed in a similar way, showing
that certain operators are vertex differential operators of higher
degree using relations between elements of homology groups, and this
can be thought of as a higher dimensional generalization of the vertex
algebra identity. At first sight there seems to be a serious problem
in carrying out this program: the higher dimensional spaces whose
homology groups we work with are not only rather complicated, but
there are far too many of them to calculate all their homology
groups. Rather surprisingly, we do not need to know the exact homology
groups, and we do not even need to find any explicit relations.  It
turns out  that the identities we want to prove
follow just from the existence of ``sufficiently many'' relations. We
will  prove the existence of enough relations  by bounding the
dimensions of the homology group and using the fact that if we have a set of
elements in the homology group with cardinality greater than its rank
then there must be a nontrivial linear relation between them. 

We let $G$ be a finite dimensional real vertex formal group.  The
underlying formal group of $G$ is the formal group of some connected
real algebraic group, and we assume that the vertex structure on $G$
is given by inverting some element of the coordinate ring of this
algebraic group (so that $K$ is of the form
$RG^*[1/p]=R[[x_1,\ldots]][1/p]$ for some polynomial $p$ in
$x_1,\ldots$).  We let $U$ stand for a unipotent algebraic subgroup of
$G$ of some dimension $n$ which is not contained in the divisor of
$p$. As $U$ is unipotent the exponential map from the Lie algebra of
$U(\C)$ to $U(\C)$ is a diffeomorphism, so that $U(\C)$ is a vector
space over $\C$ with some group structure (which is of course
different from the additive group structure if $U$ is not abelian).

The function $p$ restricts to a function on $U$ with  zero set
given by some divisor $D$.
We choose a compact $n$ cycle $C\subset U(\C)$ in the complement of
the divisor $D$, which represents an element of $H_n(U(\C)-D)$ which we
also denote by $C$.

We recall the definition of a differential operator on a ring.
If $S$ is a commutative algebra over a commutative ring $R$,
then a differential operator $D$ of order at most $n\in \Z$
is defined to be an operator which is zero if $n<0$, and if $n\ge 0$ it
is
defined to be an operator such that $[D,s]$ is a differential operator
of order at most $n-1$ for any element $s\in S$. A differential
operator $D$ is called normalized if $D(1)=0$. The differential operators
of order at most 0 are just multiplications by elements in the center of $S$,
and normalized differential operators of order at most 1
are the same as derivations of $S$, or in other words operators
such that $D(ab) = aD(b)+D(a) b $ for all $a,b\in S$. Every differential
operator can be written uniquely as the sum of
a normalized differential operator and multiplication by a constant.
If $D$ and $E$ are differential operators of orders $m$ and $n$
then $DE$ is a differential operator of order at most $m+n$, and $[D,E]$
is a differential operator of order at most $m+n-1$. In particular
the normalized differential operators of order at most 1
form a Lie algebra acting on $S$, called the Lie algebra of derivations of $S$.
A derivation of $S$ can be thought of as an infinitesimal automorphism of $S$.
We define vertex differential operators of commutative $G$ vertex algebras
in the same way. 

\proclaim Theorem 9.1.
If $a$ is an element of a commutative $G$ vertex algebra
$V$ with $G$ as above, and $n$ is the dimension of the subgroup $U$,
then $\int_C a(z) d^nz$ is a vertex differential
operator of order at most $n$. In other words if $a_0,\ldots a_n\in V$
then $$[a_0(z_0),[a_1(z_1),[\ldots,[a_n(z_n),\int_C
a(z) d^nz]\cdots]]]=0.$$

Note that $\int a(z)dz \, b(y)\ne b(y)\int a(z) dz$
even though $a(z)b(y)=b(y)a(z)$. The reason for this is that
the two integrals are taken over different cycles in the subset
where $a(z)b(y)$ is holomorphic, and these two cycles need not be
homologous if $a(z)b(y)$ has singularities.

The proof of theorem 9.1 is in two steps. We first show that the number of
linearly independent elements of the form
$$\left(\prod_{i\in S,0\le i\le N}a_i(z_i)\right)
\int_Ca(z)d^nz
\left(\prod_{i\notin S,0\le i\le N}a_i(z_i)\right)
\eqno{9.2}$$ for the $2^{N+1}$ subsets $S$ of ${0,1,2,\ldots,N}$ is
bounded by a polynomial in $N$ of degree at most $n$ (rather than the
obvious bound $2^{N+1}$) as $N$ tends to infinity. This is because the
number of linearly independent elements is bounded by the rank of
certain homology groups, and we can bound these ranks by a polynomial
in $N$.  We then show that the possible relations between these
elements are so restricted that this crude bound implies (and is even
equivalent to) a single explicit relation, which is the one given in
theorem 9.1.

\proclaim Lemma 9.3. Suppose $U(\C)$ is a complex connected
unipotent algebraic Lie group of dimension $n$
and $D$ is a closed algebraic subset. Then the rank of the homology group
$H_n(U(\C)- g_1(D)\cup g_2(D)\cup\cdots\cup g_N(D))$
for elements $g_i$ in general position
is bounded by a polynomial of degree $n$ in $N$ (whose coefficients
depend only on $U$ and $D$).

Proof. The one point compactification of $U(\C)$ is a sphere of
dimension $\dim(U(\C))=2n$, so by
Spanier-Whitehead duality it is sufficient to prove that the rank of
the homology group $H_{n-1}(g_1(D)\cup g_2(D)\cup\cdots\cup
g_N(D))$ is bounded by a polynomial of degree $n$ in $N$.  There is a
spectral sequence converging to the homology of $g_1(D)\cup
g_2(D)\cup\cdots\cup g_N(D)$ whose $E_2$ term is given by the homology
groups of all finite intersections of the $g_i(D)$'s (see [G, theorem
5.4.1 and the remarks in section 5.6]).  The rank of $H_{n-1}(g_1(D)\cup
g_2(D)\cup\cdots\cup g_N(D))$ is bounded by the sum of the ranks of
the terms $H_{n-1-(k-1)}(g_1(D)\cap g_2(D)\cap\cdots\cap g_k(D))$
of the $E_2$ term of the spectral sequence of total degree $n$, and
the rank of any homology group of any intersection of the $g_i(D)$'s
is bounded by some constant, so it is sufficient to show that the
number of ways of choosing $k$ of the $g_i$'s so that their
intersection has vanishing $H_{n-1-(k-1)}$ is bounded by a
polynomial in $N$ of degree at most $n$.  As the $g_i$'s are in
general position, the intersection of $k$ of them has dimension at
most $2(n-k)$, so that the $H_{n-k}$'th homology of this
intersection vanishes whenever $n-k>2n-2k$, or $k>n$. Hence the
number of ways of choosing $k$ of the $g_i$'s so that their
intersection has non vanishing $H_{n-1-(k-1)}$ is bounded by the number
of subsets of $N$ elements of size at most $n$, which is a polynomial
of degree $n$ in $N$.  This proves lemma 9.3.

\proclaim Lemma 9.4. The number of linearly independent
elements of the form 9.2 is bounded by a polynomial in $N$ of degree
$n$.

Proof. Each term of the form 9.2 is given by integrating
$$a(z)a_0(z_0)\cdots a_n(z_n)d^nz$$
over some $n$-cycle in $G- z_0(D)\cup\cdots\cup z_N(D)$.
By lemma 9.3 the rank of the homology group generated by
these $n$-cycles is bounded by a polynomial in $N$ of degree $n$,
which proves lemma 9.4.

In particular there must be some nonzero relation
for sufficiently large $N$, because
the number of expressions of the form 9.2 increases like $2^N$, but the number
of linearly independent ones is bounded by some polynomial in $N$.

\proclaim Lemma 9.5. Suppose that the coefficients $c_S$ are the coefficients
of a nonzero relation of smallest possible degree $N$ as above. Then
$c_S=(-1)^{|S|}c_\emptyset$ for some nonzero constant $c_\emptyset$.

Proof. We are given that
$$\sum_S c_S \left(\prod_{i\in S,0\le i\le N}a_i(z_i)\right)
\int_Ca(z)d^nz
\left(\prod_{i\notin S,0\le i\le N}a_i(z_i)\right)
=0
$$
for all $a_0,\ldots,a_N \in V$. But if we put $a_i=1$, $z_i=0$ we find
a relation of smaller degree. All coefficients of this relation must
be identically 0 as the $c_S$'s were by assumption the coefficients of
a nonzero relation of smallest degree. But the coefficients of this smaller
relation are of the form $\pm(c_S+c_{S\cup i})$ for
$S\subset\{1,\ldots,i-1,i+1,\ldots,N\}$,
so $c_S=-c_{S\cup i}$. This implies that $c_S=(-1)^{|S|}c_\emptyset$
for all $S$. This proves lemma 9.5.

In particular this shows that $\int_C a(z)d^nz$ is a differential
operator of some order. We now pin down the order more precisely
by looking more carefully at the set of all possible relations.

\proclaim Lemma 9.6. Suppose that $n'$ is the degree
of the unique nonzero relation of smallest degree, as in lemma 9.5.
Then the maximum number of linearly independent elements
of the form 9.2 is exactly equal to the number of subsets of
$\{0,1,\ldots,N\}$ of size at most $n'$ and is therefore
a polynomial in $N$ of degree exactly $n'$. 

Proof. We show that the set of elements of the form
$$\left(\prod_{i\in S,0\le i\le N}a_i(z_i)\right)
\int_Ca(z)d^nz
\left(\prod_{i\notin S,0\le i\le N}a_i(z_i)\right)
$$ with $|S|\le n'$ form a maximal linearly independent set of functions
of the $a_i(z_i)$'s. In fact by using the relation of lemma 9.5 we see
that any relation can be written as a sum of relations with at most
$n'$ factors to the left of the integral. To complete the proof of the
lemma we have to show there are no nontrivial linear relations between
these terms. Suppose there is a nontrivial relation with coefficients
$c_S$, with $c_S=0$ if $|S|> n'$. We can assume that we
have chosen this relation so that the maximum value of
$|S|$ with $c_S\ne 0$ is as small as possible. Now set
$a_i=1$, $z_i=0$ for $i\notin S$. Then we get a relation
of degree at most $n'$ with a nonvanishing coefficient of
$\left(\prod_{i\in S,0\le i\le N}a_i(z_i)\right)
\int_Ca(z)d^nz$, which is impossible.
Hence the maximal number of linearly independent elements
is exactly equal to the number of subsets of size at most $n'$
of a set of size $N$, which is a polynomial in $N$ of degree
exactly $n'$. This proves lemma 9.6.

We can now complete the proof of theorem 9.1. By lemma 9.4
the maximal number of linearly independent elements of the form 9.2 is
bounded by a polynomial of degree at most $n$. In particular there
must exist some relation of minimal degree $n'$ by lemma 9.5.
By lemma 9.6 we see that $n'\le n$, so by lemma 9.5 the relation 9.1
holds, in other words $\int_Ca(z)d^nz$ is a vertex differential operator
of order at most $n'\le n$. This proves theorem 9.1.

Theorem 9.1 can easily be generalized in several ways as follows.
\item{1} The
condition that $C$ should lie inside some unipotent subgroup $U$ can
be removed; it is put in only because it simplifies the proof slightly
and is satisfied in all the examples we use later.
\item{2} The restriction
that $G$ is finite dimensional is also usually unnecessary; for example
we could take $G$ to be the formal group of the Virasoro algebra and
take $U$ to correspond to the 1 dimensional group generated by
$L_{-1}$.
\item{3} We can also look at vertex differential operators $a(x)$
and find that $\int_C a(z)d^nz$ is a vertex differential operator of
order at most the order of $a$ plus the dimension of $U$.
\item{4} We can also
look at vertex differential operators $a(x,y,z,\ldots) $ in several
parameters, and find that
$b(y,z,\ldots)=\int_Ca(x,y,z,\ldots)d^nx$ is a vertex
differential operator with order as above.
\item{5} We can include a singular function $f$ of several variables
in the integrand without affecting the argument. For example
if $a$ is an element of a commutative $G$ vertex algebra
$V$ with $G$ as above, and $n$ is the dimension of the subgroup $U$,
then $\int_C a(z)f(z,z_0,\ldots,z_n) d^nz$ is a vertex differential
operator of order at most $n$.
\item{6} We do not need to restrict
the vertex structure on $G$ to be given by inverting a function
$p$ on an algebraic group, and we can allow things like the
vertex structure generated by $f^{1/n}$ for some integer $n$.
\item{7} Finally the restriction to vertex  groups over the reals is
unnecessary and  theorem 9.1 can be
generalized to vertex  groups over any field
$R$ by using \'etale cohomology instead of singular cohomology.

\proclaim  10.~$G$ Vertex algebras and the Yang-Baxter equation.

In this section we show how to construct examples of 
$G$ vertex algebras from solutions of the Yang-Baxter equation.
The idea is to start with an  $G$ algebra,
and deform it using a solution of the Yang-Baxter equation into
a commutative $G$ vertex algebra. (Note that in many examples
we end up with a commutative $G$ vertex algebra, even though
the algebra we start with is not commutative!)

We first consider the case of associative algebras $V$.
We say that a linear map $R$ from $V\otimes V$ to $V\otimes V$
is an $R$ matrix if it satisfies the following conditions:
\item{1.} $R_{12}R_{13}R_{23}=R_{23}R_{13}R_{12} $
(Yang-Baxter equation.) Both sides act on $V\otimes V\otimes V$,
and $R_{ij}$ means $R$ acting on the $i$'th and $j$'th factors
of the tensor product.
\item{2.} $R(1\otimes v)=1\otimes v$, $R(v\otimes 1)=v\otimes 1$.
\item{3.} $R_{12}m_{12}=m_{12} R_{23}R_{13}$,
$R_{12}m_{23}=m_{23} R_{12}R_{13}$, where $m_{ij}$ is the product 
map from $V\otimes V\otimes \cdots V$ to $V\otimes \cdots V$ given by
multiplying the $i$'th and $j$'th factors. 

Note that we have added extra conditions saying that the $R$ matrix
is ``compatible'' with the algebra structure.

\proclaim Lemma 10.1. Suppose $V$ is a commutative algebra
and $R$ is an $R$ matrix for $V$. Then the bilinear
map $m_{12}R_{12}$ is another associative algebra structure on
$V$ (with the same identity element 1).

Proof. The only nontrivial thing to check is associativity,
which follows from
$$\eqalign{
&m_{12}R_{12}m_{23}R_{23}\cr
=&m_{12}m_{23}R_{12}R_{13}R_{23}\cr
=&m_{12}m_{12}R_{23}R_{13}R_{12}\cr
=&m_{12}R_{12}m_{12}R_{12}.\cr
}$$

Now we assume that $V$ is acted on by a group $G$. The group
$G\times G$ then acts on $V\otimes V$ and on $Hom(V\otimes V,V\otimes V)$.
We denote the image of the matrix $R$ under $g\times g'$ by $R^{g,g'}$.
We assume that $R$ satisfies the following axioms. The first three
are the obvious analogues of the ones above, the fourth is just
the definition of the action of $G\times G$ on $R$, and the fifth is
a sufficient condition for the $G$-invariance of
the new algebra structure on $V$.

\item{1.} $R_{12}^{g_1,g_2}R_{13}^{g_1,g_3}R_{23}^{g_2,g_3}
=R_{23}^{g_2,g_3}R_{13}^{g_1,g_3}R_{12}^{g_1,g_2} $
(Yang-Baxter equation.)
\item{2.} $R^{g_1,g_2}(1\otimes v)=1\otimes v$, $R^{g_1,g_2}(v\otimes
1)=v\otimes 1$.
\item{3.} $R_{12}^{g_1,g_2}m_{12}=m_{12} R_{23}^{g_1,g_2}R_{13}^{g_1,g_2}$,
$R_{12}^{g_1,g_2}m_{23}=m_{23} R_{12}^{g_1,g_2}R_{13}^{g_1,g_2}$.
\item{4.} $(R^{g_1,g_2}(u\otimes v))^{g_1',g_2'}
=R^{g_1g_1',g_2g_2'}((u\otimes v)^{g_1',g_2'})$
\item{5.} $R^{gg_1,gg_2}=R^{g_1,g_2}$

If we want we can use the fifth axiom to define
$R^g=R^{1,g}$ so that $R^{g_1,g_2}=R^{g_1^{-1}g_2}$, so
$R$ really only depends on one element of $G$.

As before, we can use $R^{g_1,g_2}$ to twist a $G$ invariant algebra
product on $V$ to get a new $G$-invariant algebra product.

We can now try to construct $G$ vertex algebras in the same way,
except that the matrix $R^{g_1,g_2}$ is a singular function
of $g_1,g_2$; more precisely, $R^g$ should be an element
of $Hom(V\otimes V, Hom(H, V\otimes V)\otimes_{H^*}K)$.
The quantum groups literature contains many examples
of solutions of equation 1 above with singularities for $g_1,g_2$
complex numbers, or more generally elements of some
Riemann surface. These do not usually satisfy equation
3 because the space $V$ is usually taken to be a finite dimensional
vector space rather than an algebra, but equation 3 shows that there
is at most one extension of $R$ to any algebra generated
by this finite dimensional space, such as the tensor algebra.
This gives a large number of $G$ vertex algebras
constructed using solutions of the Yang-Baxter equation.

{\bf Example 10.2.} We take $G$ as in example 3.4, and show how to construct
the vertex
algebra of a lattice $L$ using the method above. For simplicity
we assume all inner products in $L$ are even; the general case can be done
in a similar way using a twisted group ring of $L$.
We take $V$ to be the universal commutative $H$ algebra generated by the group
ring $\Z[L]$ of $L$. There is a singular
solution to the Yang-Baxter equation on $\Z[L]$,
given by $$R^{x,y}(e^a\otimes e^b)=(x-y)^{(a,b)}e^a\otimes e^b.$$
This extends uniquely to an $R$ matrix on $V$ satisfying the
conditions above, which can be used to make $V$ into a vertex algebra.
This is the usual vertex algebra of an even lattice.

\proclaim 11.~Cohomology of $G$ vertex algebras.

We have seen earlier that many $G$ vertex algebras can be constructed
by deforming the product on a $G$-algebra. In this section we will show
how to describe the infinitesimal deformations
of a $G$-algebra into a $G$ vertex algebra using
$G$-equivariant cohomology of associative algebras.

We start by recalling how to classify
the infinitesimal deformations of an associative algebra $V$
over a commutative ring $R$. This means that
we want to define an associative algebra structure
on $V[\epsilon]$ over $R[\epsilon]$, where $\epsilon^2=0$.
For simplicity we will work with associative algebras
without identity elements; the main difference below if we
use algebras with identity elements is that we should use
normalized cochains instead of cochains.
If we write the new algebra product as
$ab+\epsilon f(a,b)$ then associativity is equivalent to
$$af(b,c)-f(ab,c)+f(a,bc)-f(a,b)c=0$$
so $f$ is just a 2-cycle for the standard complex
used to calculate the associative algebra cohomology $H^2(V,V)$.
(Recall that the $n$-cochains of the standard resolution
for calculating the Hochschild cohomology groups
$H^n(V,M)$ for a 2 sided module $M$ over
the associative algebra $V$ are the multilinear maps
from $A,A,\ldots, A$ to $M$; see [C-E].)
The ``trivial'' deformations are those that are induced
by an infinitesimal automorphism $a\mapsto a+\epsilon g(a)$
of the underlying $R$-module, when the corresponding
2-cocycle $f$ is given by $f(a,b)=ag(b)-g(ab)+g(a)b$,
in other words $f$ is just the coboundary of the 1-cochain $g$.
Therefore the infinitesimal deformations of the algebra
$V$ are classified by the cohomology group $H^2(V,V)$
of the associative algebra $V$ with coefficients in the 2-sided module $V$.

Now suppose that $V$ is acted on by a group (or cocommutative Hopf algebra)
$G$. Then we can calculate the infinitesimal deformations
of the $G$-algebra $V$ by using $G$-equivariant cohomology
groups $H^2_G(V,V)$. The $G$-equivariant cohomology
groups are calculated (and defined) by replacing the (standard) cochains
above by the $G$-invariant cochains.

Finally suppose that $G$ is a vertex group.
We define the cohomology groups $H^n_G(V,M)$
by using the $G$-invariant singular multilinear functions
from $V,V,\ldots V$ to $M$ in place of the multilinear
functions from $V,V\ldots,V$ to $M$. (The boundary operator
is defined by the same formula as for associative algebras.)
We then find by the argument above that
the is a map from $H^2_G(V,V)$ for an associative $G$-algebra $V$
to $G$ vertex algebras over $R[\epsilon]$. The group $H^2_G(V,V)$
can be thought of a roughly the tangent space at $A$
of the moduli space of $G$ vertex algebra structure on $A$.
Some examples of nontrivial elements of $H^2_G(V,V)$
can easily be obtained from the $G$ vertex algebras
of generalized free fields, because we can just take
a 2-point function $\Delta$ of the form $\epsilon/(x-y)^2$
for example.

We can also define other sorts of cohomology groups involving vertex
groups $G$, such as vertex group cohomology.  We first note that most
definitions in group cohomology (in particular the homogeneous and
inhomogeneous standard complexes) can easily be extended to the case
of arbitrary cocommutative Hopf algebras.  We can now define the
vertex group cohomology groups $H^n(G,M)$ for a vertex group $G$ and a
$G$-module $M$ by replacing the spaces of multilinear maps
$Hom(H,H\ldots, H, M)$ in the homogeneous standard complex (where $H$
is the underlying cocommutative Hopf algebra) by the corresponding
spaces of singular multilinear maps.

\proclaim 12.~Relation to quantum field theory.

There are several mathematical structures closely
related to quantum field theory;  for example, Wightman's axioms [S-W]
and some closely related variations [G-J], Segal's axioms
for a topological field theory [Se], and $G$ vertex algebras
as described above. We give an informal description of
the relation between $G$ vertex algebras and the other theories above.
The four theories above are closely related but not exactly the same.

The relation with Wightman's axioms is easiest to describe.
A quantum field theory satisfying Wightman's axioms is determined
by its correlation functions [S-W], and if these correlation functions
have ``good'' operator product expansions then they
are the correlation functions of some $G$ vertex algebra,
with $G$ the vertex group whose underlying Hopf algebra
is the universal enveloping algebra of the Poincar\'e algebra, 
and where we allow some sort of singularities on the light cone.
Unfortunately it is
difficult to decide when a $G$ vertex algebra comes from
a quantum field theory satisfying the Wightman axioms. The main problem is
in reconstructing the Hilbert space. It is usually not too hard to
construct a real vector space with a symmetric bilinear form on it,
but deciding when this form is positive (semi) definite is
usually rather hard. (For example a special case of this is the
problem of deciding which highest weight representations of the Virasoro
algebra are unitary.)

Note that expressions like $\phi_1(x_1)\phi_2(x_2)$ are interpreted in
quite different ways in quantum field theories and $G$ vertex
algebras: in quantum field theories we think of the $\phi_i$'s as
distribution valued operators defined on a manifold containing points
$x_i$, and the product is a product of operators.  In $G$ vertex
algebras the elements $x_i$ are thought of as elements of some group
and $\phi(x)$ should be thought of as the transform of the element
$\phi$ under the action of the group element $x$. The product does not
always make sense for fixed $x_1$ and $x_2$, but is only defined as
some sort of singular function of $x_1$ and $x_2$.

The relation with quantum field theories that occur in physics is similar:
provided that good operator product expansions exist, there
is often a $G$ vertex algebra with the same correlation functions.
The main problem is that most realistic quantum field theories are only
defined at the level of perturbation theory, in other words,
the correlation functions are formal power series in the coupling constants
that (probably) do not converge for any nonzero values.
We can get round this by the following trick: we define the $G$ vertex algebra
over the ring of formal power series in the coupling constants.
So many quantum field theories are now well defined mathematical objects:
they are $G$ vertex algebras over formal power series rings.
Notice that this does not solve the important problem of making sense
of these theories non perturbatively; all we have done is
change the question of existence of quantum field theories
into an equally hard question about properties of
$G$ vertex algebras over formal power series rings.
More precisely, we would really like to construct some sort
of moduli space of $G$ vertex algebras, such that at points of its
compactification the $G$ vertex algebras somehow degenerate into the
formal power series $G$ vertex algebras above.

The relation of $G$ vertex algebras with Segal's topological field theories is
harder to describe; in fact, there seems to be no particularly
easy way to go between them. The reason for this seems
to be that $G$ vertex algebras like to work with groups $G$
and work best when spacetime can be regarded as a group (for example,
if spacetime is flat). On the other hand, Segal's axioms work
with arbitrary manifolds, most of which have nothing to do with groups.
A good example is given by 2 dimensional conformal field theories.
It is well known that vertex operator algebras (which are more
or less vertex $Virasoro$-algebras) are good at describing the
genus 0 part of a conformal field theory, but are  bad at describing
the higher genus part (although the genus 1 case has been pushed through
by Zhu [Z]). This is because a genus 0 Riemann surface is
more or less the additive group $\C$ (at least if a point is missed out),
while Riemann surfaces of genus greater than 1 are not directly related to
groups.
The additive group $\C$ is of course more or less the ``group'' that
acts on vertex algebras.

To summarize,  quantum field theories satisfying the Wightman axioms,
topological field
theories, and $G$ vertex algebras are related but different
mathematical structures, each of which captures part but not all of
the notion of a quantum field theory. Algebraic quantum field theories
emphasize Lorentz invariance, locality, and unitarity, but have the
disadvantage that it is extraordinarily difficult to construct
nontrivial examples of them.  $G$ vertex algebras emphasize locality,
Poincare invariance, and the operator product expansion,
and can cope with the Feynman path integral (regarded
as a trace on the $G$ vertex algebra), but are not
very good at dealing with unitarity or curved spacetimes. Segal's axioms
emphasize the Feynman path integral and are particularly good
at dealing with curved spacetimes.

\proclaim References.

\item{[A]} E. Abe, Hopf algebras, Cambridge University Press, 1980.
\item{[B-D]} A. Beilinson, V Drinfeld, Chiral algebras I. 1996(?) preprint.
\item{[B]}{R. E. Borcherds,  Vertex algebras, Kac-Moody algebras,
   and  the monster. Proc. Natl. Acad. Sci.  USA. Vol. 83 (1986) 3068--3071.}
\item{[C-E]} H. Cartan, S. Eilenberg, Homological algebra,
Princeton U. P. 1956.
\item{[G-J]} J. Glimm, A. Jaffe, Quantum physics, Springer Verlag 1981.
\item{[G]} R. Godement, Th\'eorie des faisceaux, Hermann 1964.
\item{[K]} V. G. Kac, Vertex algebras for beginners,
University lecture series vol. 10, A.M.S. 1996.
\item{[L-Z]} B. Lian, G. Zuckerman, New perspectives on the
BRST-algebraic structure in string theory, hep-th/9211072,
Communications in Mathematical Physics 154, (1993) 613--64, and
Moonshine cohomology, q-alg/950101, Finite groups and Vertex Operator
Algebras, RIMS publication (1995) 87-11.
\item{[Se]} Segal, Graeme,  Geometric aspects of quantum field theory.
Proceedings of the International Congress
of Mathematicians, Vol. I, II (Kyoto, 1990),
1387--1396, Math. Soc. Japan, Tokyo, 1991.
\item{[S]} J-P. Serre,  ``Lie algebras and Lie groups'', 1964 lectures given at
Harvard University. Second edition. Lecture Notes in Mathematics, 1500.
Springer-Verlag, Berlin, 1992. 
\item{[So]} Y. Soibelman, Meromorphic tensor categories,
preprint q-alg/9709030.
\item{[S-W]} R. F. Streater, A. S. Wightman, PCT, Spin, statistics,
and all that. Addison-Wesley 1964.
\item{[Z]} Zhu, Yongchang,
Modular invariance of characters of vertex operator algebras.
J. Amer. Math. Soc. 9
(1996), no. 1, 237--302.
\bye